\documentclass[11pt,a4paper,twoside]{article}
\usepackage{amsfonts}
\usepackage{amsthm,amsmath,amssymb}
\usepackage{color}
\usepackage{cite}
\def\trait #1 #2 #3 {\vrule width #1pt height #2pt depth #3pt}
\newtheorem{teor}{Theorem}[section]
\newtheorem{defin}[teor]{Definition}
\newtheorem{lemm}[teor]{Lemma}
\newtheorem{osse}[teor]{Remark}
\newtheorem{prop}[teor]{Proposition}
\newtheorem{defi}[teor]{Definition}
\newtheorem{coro}[teor]{Corollary}
\newtheorem{prob}[teor]{Problem}
\newtheorem{hypo}[teor]{Hypothesis}
\newcommand{\bele}{\begin{lemm}\begin{sl}}
\newcommand{\enle}{\end{sl}\end{lemm}}
\newcommand{\bedef}{\begin{defi}\begin{sl}}
\newcommand{\eddef}{\end{sl}\end{defi}}
\newcommand{\bete}{\begin{teor}\begin{sl}}
\newcommand{\ente}{\end{sl}\end{teor}}
\newcommand{\beos}{\begin{osse}\begin{rm}}
\newcommand{\eddos}{\end{rm}\end{osse}}
\newcommand{\bepr}{\begin{prop}\begin{sl}}
\newcommand{\empr}{\end{sl}\end{prop}}
\newcommand{\bepro}{\begin{prob}\begin{rm}}
\newcommand{\empro}{\end{rm}\end{prob}}
\newcommand{\bede}{\begin{defin}\begin{sl}}
\newcommand{\edde}{\end{sl}\end{defin}}
\newcommand{\beco}{\begin{coro}\begin{sl}}
\newcommand{\enco}{\end{sl}\end{coro}}
\newcommand{\behy}{\begin{hypo}\begin{sl}}
\newcommand{\enhy}{\end{sl}\end{hypo}}

\newcommand{\RR}{\mathbb{R}}

\newcommand{\beeq}[1]{\begin{equation}\label{#1}}
\newcommand{\eddeq}{\end{equation}}
\newcommand{\beeqa}[1]{\begin{eqnarray}\label{#1}}
\newcommand{\eddeqa}{\end{eqnarray}}
\newcommand{\beal}[1]{\begin{align}\label{#1}}
\newcommand{\eddal}{\end{align}}
\newcommand{\bespl}[1]{\begin{split}\label{#1}}
\newcommand{\edspl}{\end{split}}
\newcommand{\bega}[1]{\begin{gather}\label{#1}}
\newcommand{\edga}{\end{gather}}
\newcommand{\grad}{\nabla}
\newcommand{\beeqax}{\begin{eqnarray*}}
\newcommand{\eddeqax}{\end{eqnarray*}}
\newcommand{\curlu}{\nabla \times \mathbf{u}}
\numberwithin{equation}{section}
\def\qed{\ifmmode   \else \leavevmode\unskip\penalty9999 \hbox{}\nobreak\hfill
  \fi
  \quad\hbox{\hskip.5em\vrule width.4em height.6em depth.05em\hskip.1em}}
\def\endproofsym{\qed}

\renewenvironment{proof}[1][Proof]{\trivlist\item[\hskip\labelsep{\hskip0pt
        {\normalfont\scshape#1.}\hskip .321429\parindent}]\ignorespaces}
{\endproofsym\endtrivlist}
\def\endnobox{\def\endproofsym{}\end{proof}\def\endproofsym{\qed}}

\newcommand{\no}{\nonumber}
\newcommand{\beeqao}{\begin{eqnarray}\no}
\newcommand{\bealo}{\begin{align}\no}
\newcommand{\besplo}{\begin{split}\no}
\newcommand{\begao}{\begin{gather}\no}

\newcommand{\vc}[1]{{\boldsymbol #1}}
\newcommand{\dt}{\partial_t}

\newcommand{\io}{\int_\Omega}

\newcommand{\bsig}{\boldsymbol{\sigma}}

\newcommand{\bd}{\mathbf{d}}

\newcommand{\Grad}{\nabla}

\newcommand{\ub}{{\mathbf{u}}}
\newcommand{\vb}{{\mathbf{v}}}

\def\fine{\hfill\kern4pt \vrule height4pt depth0pt width4pt }

\def\dive{\mbox{\rm div\,}}

scaled \magstep1 

\allowdisplaybreaks[4]

\begin{document}

\title{Global weak solution and blow-up criterion of the general Ericksen--Leslie
system for nematic liquid crystal flows}

\author{Cecilia Cavaterra
\\
Dipartimento di Matematica, Universit\`{a} degli Studi di Milano\\
Via Saldini 50, 20133 Milano, Italy
\\
cecilia.cavaterra@unimi.it
\and
Elisabetta Rocca
 \thanks{The author was supported by the FP7-IDEAS-ERC-StG
Grant \#256872 (EntroPhase)}\\
Dipartimento di Matematica, Universit\`{a} degli Studi di Milano\\
Via Saldini 50, 20133 Milano, Italy
\\
elisabetta.rocca@unimi.it
\and
Hao Wu
 \thanks{Corresponding author. The author was partially supported by NSF of China 11001058, Specialized Research Fund for
the Doctoral Program of Higher Education and "Chen Guang" project supported by Shanghai Municipal Education
Commission and Shanghai Education Development
Foundation. Part of the work was done during his staying in Milan in May 2012 within the FP7-IDEAS-ERC-StG
Grant \#256872}\\
School of Mathematical Sciences\\
Shanghai Key Laboratory for Contemporary Applied Mathematics\\
Fudan University\\
Han Dan Road 220, 200433 Shanghai, China\\
haowufd@yahoo.com
}
\date{\today}
\maketitle

\begin{abstract}
\noindent
In this paper we investigate the three dimensional general Ericksen--Leslie (E--L)
system with Ginzburg-Landau type approximation modeling nematic liquid crystal flows.
First, by overcoming the difficulties from lack
of maximum principle for the director equation and high order nonlinearities for the stress tensor, we prove existence of global-in-time weak solutions under physically
meaningful boundary conditions  and suitable assumptions on the Leslie coefficients, which ensures that the total
energy of the E--L system is dissipated.
Moreover, for the E--L system with periodic boundary conditions, we prove the local well-posedness of classical solutions under the so-called Parodi's relation and establish a
blow-up criterion in terms of the temporal integral of both the maximum norm of the curl of the velocity field and the maximum norm of the gradient of the liquid crystal director field.

\smallskip


\smallskip

\noindent
{\bf AMS (MOS) subject clas\-si\-fi\-ca\-tion:}~~35B44, 35D30, 35K45, 35Q30, 76A15.
\end{abstract}


\section{Introduction}
\label{sec:intro}

The hydrodynamic theory of liquid crystals
due to Ericksen and Leslie was developed around 1960's
\cite{E61,E62,Le68}.
The Ericksen--Leslie system for liquid crystal flows was derived from the macroscopic point of view and was
very successful in understanding the coupling between the
velocity field and the molecular director field, especially in liquid crystals of nematic type.
Assume that the material occupies a bounded spatial domain $\Omega\subset \RR^3$ with smooth boundary $\Gamma$.
We denote by  $\ub=(u_1,u_2,u_3)^T$ the velocity field of the flow and by $\bd=(d_1,d_2,d_3)^T$ the director field,
which stands for the averaged macroscopic/continuum molecular orientation in $\RR^3$.
The resulting PDE system that we are going to consider in $Q_T :=(0,T)\times \Omega$
can be written as \cite{LL01,WXL12}:
 \begin{eqnarray}
 \ub_t+(\ub\cdot\nabla) \ub+\nabla P&=&-\nabla\cdot(\nabla \bd
\odot\nabla \bd) +\nabla\cdot\bsig,
 \label{e1}\\
\nabla\cdot \ub&=&0, \label{e2}\\
\bd_t+(\ub\cdot\nabla )\bd-\omega
\bd+\frac{\lambda_2}{\lambda_1}A\bd&=&-\frac{1}{\lambda_1}\left(\Delta
\bd-\Grad_{\bd} W(\bd)\right). \label{e3}
\end{eqnarray}
Equations \eqref{e1} and \eqref{e3} represent the
conservation of linear momentum and angular momentum, respectively. Here, we consider the
flow of an incompressible material that satisfies the incompressiblility condition \eqref{e2}.
$P$ is a scalar function representing the hydrodynamic pressure (including the hydrostatic part and the induced
elastic part from the director field). The term $-\Grad \bd \odot \Grad \bd$ is the elastic (Ericksen)
stress tensor, where $\Grad \bd \odot \Grad \bd$ denotes the $3\times 3$ matrix whose $(i,j)$-th entry is given by
$\nabla_i \bd \cdot \nabla_j \bd$, for $1\leq i,j\leq 3$. The function $W$ in \eqref{e3} penalizes the deviation of the length $|\bd|$ from the value 1,
which is due to liquid crystal molecules being of similar size (cf. \cite{LinLiusimply,LL01,WZZ12}).
The notation $\Grad_{\bd}$ represents the gradient with respect to the variable $\bd$.
A typical example of $W$ is the Ginzburg--Landau approximation that has a double well structure \cite{LinLiusimply, LL01}:
\begin{equation}
W(\bd)=\frac{1}{4\varepsilon^2}(|\bd|^2-1)^2.\label{GL}
\end{equation}
 The symbol $\Grad^T$ indicates the transpose of the gradient and
$$
 A=\frac12(\nabla \ub+\nabla^{T}\ub),\
\,\; \omega=\frac12(\nabla \ub-\nabla^{T}\ub)$$
represent the rate of strain tensor  and the skew-symmetric part of the strain rate, respectively.
For the sake of simplicity, we denote by
\begin{equation}
\dot{\bd}=\bd_t+( \ub\cdot\nabla)\bd, \ \,\; \mathcal{N}=\dot{\bd}-\omega\,\bd
= -\frac{\lambda_2}{\lambda_1}A\bd -\frac{1}{\lambda_1}\left(\Delta \bd-\Grad_{\bd} W(\bd)\right)
\label{def1}
\end{equation}
the material derivative of $\bd$ (transport of center of mass) and the rigid
rotation part of the changing rate of the director by fluid vorticity.
The kinematic transport $\lambda_1 \mathcal{N}+\lambda_2 A\bd $ for the molecule
director represents the effect of the macroscopic flow field on the microscopic structure.
The material coefficients $\lambda_1$ and $\lambda_2$ reflect the molecular shape (Jeffreys orbit) and how slipper
the particles are in the fluid  (cf. \cite{Jeff,sunliu}).
The viscous (Leslie) stress tensor $\bsig$ has the following form \cite{dG,Le79}:
\begin{equation}
\bsig=\mu_1(\bd^TA \bd)\bd\otimes \bd+\mu_2\mathcal{N}\otimes \bd+\mu_3\bd\otimes \mathcal{N}+
\mu_4A+\mu_5(A\bd)\otimes \bd+\mu_6\bd\otimes(A\bd), \label{v5}
\end{equation}
where $\otimes$ stands for the usual Kronecker product, i.e.,
$(\mathbf{a}\otimes\mathbf{b})_{ij}:=a_ib_j$, for $1\leq i,j\leq 3$.
The six independent coefficients $\mu_1,...,\mu_6$, which may depend
on material and temperature, are called Leslie coefficients. These
coefficients are related to certain local correlations in the fluid (cf. e.g., \cite{dG} for more details).

The system \eqref{e1}--\eqref{e3} should be supplemented with suitable initial and boundary conditions.
Here, we assume that it is subject to the following  initial conditions
 \begin{equation} \ub|_{t=0}=\ub_0 \ \ \mbox{with} \ \ \nabla\cdot
\ub_0=0, \quad \quad \bd|_{t=0}=\bd_0, \ \ \mbox{in} \ \Omega,
 \label{IC}
 \end{equation}
and the homogeneous Dirichlet
boundary condition for the velocity field
\begin{equation}\label{slip}
\ub  = \mathbf{0} , \quad  \mbox{on} \ (0,T)\times\Gamma,
\end{equation}
together with the nonhomogeneous Dirichlet
boundary condition for $\bd$  (constant in time for simplicity)
\begin{equation}\label{diri}
\bd|_{\Gamma} = \bd_0|_{\Gamma},  \quad \mbox{on} \ (0,T)\times\Gamma.
\end{equation}

The reformulated E--L system \eqref{e1}--\eqref{e3} with Ginzburg--Landau approximation $W$ was introduced in \cite{LL01}
motivated by the work on the harmonic heat flow (see, e.g., \cite{LinLiusimply,S02} for a detailed analysis on a simplest
liquid crystal system, in which the Leslie stress $\bsig$ is neglected except the fluid viscosity term involving $\mu_4$).
It can also be derived via an energetic variational approach (cf. \cite{WXL12}) based on the basic energy law of the system
under proper assumptions on the Leslie coefficients (e.g., the Parodi's relation, see \eqref{lam2} below).
Due to the complicated mathematical structure of the E--L system, most of the previous works were restricted to certain simplified
versions, see for instance, \cite{CR, FR, GrasWu, LinLiusimply,LL01, LLZ07, LS, LZC12, prslong, S02, sunliu, wuxuliu} and the references therein.
In particular, existence of weak solutions to the E--L system \eqref{e1}--\eqref{e3} (with Ginzburg--Landau approximation) subject
to boundary conditions of Dirichlet type \eqref{slip}--\eqref{diri} was obtained in \cite{LL01} under the specific assumption that $\lambda_2=0$.
Physically this assumption indicates that the stretching effect due to the flow field is neglected, which is more feasible for small
molecules  \cite{LinLiusimply}. From the mathematical point of view, it brings great convenience since a weak maximum principle for
$|\bd|$ holds (cf. \cite[Theorem 3.1]{LL01}, also \cite{LinLiusimply}).
 When the stretching effect is taken into account, several results have been obtained. For instance, by taking
   \begin{eqnarray}
 &&  \lambda_1=-1,\quad \lambda_2=2\alpha-1,\no\\
&& \mu_1=0,\ \ \mu_2=-\alpha, \ \ \mu_3=1-\alpha,\ \mu_4=\nu,\no\\
&&   \mu_5=\alpha(2\alpha-1),\ \ \mu_6=(\alpha-1)(2\alpha-1), \no
 \end{eqnarray}
for some constants $\alpha\in [0,1]$ and $\nu>0$ (cf. \cite[Remark 3.2]{WXL12}), we can formulate a simplified version of the E--L system
\eqref{e1}--\eqref{e3} such that the molecule has a general ellipsoid shape (cf. e.g., \cite{LLZ07,sunliu}).
The parameter $\alpha$ is related to the shape of the liquid crystal molecules. For instance, the spherical, rod-like and disc-like liquid crystal
molecules correspond to the cases $\alpha= \frac{1}{2}, \, 1$ and $0$, respectively (cf. \cite{Jeff,sunliu}).
Under the above choice of coefficients, if $\alpha\neq \frac12$ then $\lambda_2\neq 0$ and the stretching effect is kept into account.
We refer to \cite{sunliu} for the case with rod-like molecule $(\alpha=1)$, where existence of solutions in the regularity class
 \begin{equation}
 \ub \in L^\infty(0,T; H^1)\cap L^2(0,T; H^2), \quad \bd\in L^\infty(0,T; H^2)\cap L^2(0,T; H^3)
 \label{morereg}
 \end{equation}
was obtained in 2D, or in 3D under the assumption that the viscosity coefficient $\nu$ is sufficiently large with respect to proper Sobolev norms
of the initial data. See also \cite{ZF10} for certain blow-up criteria of the local classical solutions. For the general case $\alpha\in [0,1]$, we refer to
\cite{GrasWu,prslong,wuxuliu} for detailed analysis on well-posedness as well as long-time behavior of the system.
As was pointed out in \cite{sunliu}, due to the lack of maximum principle for the director equation when $\lambda_2\neq 0$, it seems to be difficult
to define (energy bounded) weak solutions in the usual sense as in \cite{LinLiusimply,LL01} and one has to deal with more regular solutions like \eqref{morereg}.
Recently,  existence of certain suitably defined \emph{weak solutions} with energy bound was obtain in \cite{CR}, with the help of an appropriate choice of the
test functions that leads to a rigorous weak formulation of the problem. This method has also been used in \cite{ffrs},
where the authors prove the existence of weak solutions for a non-isothermal variant of the liquid crystal system in \cite{CR, prslong, LLZ07, sunliu, wuxuliu}
with Neumann (for $\bd$) and complete slip (for $\ub$) boundary conditions.

We note that the liquid crystal system studied in \cite{CR,prslong,sunliu,wuxuliu} is much simpler than the general E--L system \eqref{e1}--\eqref{e3}, because
the six independent Leslie coefficients $\mu_1, ..., \mu_6$ are only characterized by two  parameters $\alpha, \nu$.
However, it has been shown that relations between those coefficients are crucial to retain certain basic properties of the general E--L system \eqref{e1}--\eqref{e3}
(with penalty), e.g., the dissipative basic  energy law \cite{LL01, WXL12}.
In particular, when the stretching effect is taken into account (i.e., $\lambda_2$ is not necessary to be zero), in \cite{WXL12} the authors provided sufficient
conditions on the Leslie coefficients to ensure that the total energy of the E--L system \eqref{e1}--\eqref{e3} is dissipated  (cf. Section 2.1) and
they discussed the important role of the so-called Parodi's relation (cf. \eqref{lam2}) in the well-posedness and stability of the system.
Besides, well-posedness and long-time behavior of E--L system \eqref{e1}--\eqref{e3} were obtained under suitable assumptions on the Leslie coefficients
(for instance, the fluid viscosity $\mu_4$ is sufficiently large).
They work within the same regularity class (cf. \eqref{morereg}) for the solution $(\ub, \bd)$ as in \cite{GrasWu,sunliu,wuxuliu},
due to the same technical difficulty, i.e, the lack of control on the $L^\infty$-norm of the director field $\bd$.
It is worth mentioning that in the recent work \cite{WZZ12} a sufficient and necessary condition on the Leslie coefficients was shown to ensure that the energy
of the original E--L system with constraint $|\bd|=1$ is dissipated,
and local well-posedness as well as global well-posedness for small initial data were obtained. We refer to the recent works \cite{CLL, HW12, HWW12, LinLinWang10, Wang11, WD11, XZ12}
and the references cited therein for mathematical results on various (simplified) liquid crystal systems under the constraint $|\bd|=1$.

The first aim of this paper is to show that, when the stretching effect is taken into account (namely, $\lambda_2\neq  0$ is allowed) and without any restriction on
the size of the fluid viscosity $\mu_4$ as well as the initial data (as it is for the uncoupled 3D Navier--Stokes system), it is possible
to obtain the existence of certain suitably defined \emph{weak solutions} with finite energy for the E--L system \eqref{e1}--\eqref{e3} (see Theorem \ref{theo1}).
The key point relies on an appropriate choice of the test functions leading to a rigorous weak formulation of the system (cf. Definition \ref{defweak}).
This technique is analogous to the one first used in \cite{CR} for the simplified system, but with more delicated estimates. In order to treat the nonlinear
stretching terms as well as  higher-order stress terms in the system,  we make use of the a double-level approximate scheme that is a combination of a suitable
Faedo--Galerkin approximation and a regularization procedure.
In particular, a nonstandard but physically meaningful regularization of the momentum equation is introduced by adding to it a $r$-Laplacian operator acting on
the velocity $\ub$, i.e., we add in the stress tensor a term of the type $|\Grad \ub|^{r-2}\Grad \ub$. Then we shall pass to the limit after obtaining suitable
uniform estimates only from the basic energy law.

In the second part of the paper, we will consider local classical solutions to the E--L system \eqref{e1}--\eqref{e3}. Under the assumption that the Parodi's
relation (cf. \eqref{lam2}) is satisfied (namely, under the set of assumptions \textbf{Case 1} below), we first prove the existence and uniqueness of local
classical solutions to the E--L system \eqref{e1}--\eqref{e3} subject to periodic boundary conditions (cf. Theorem \ref{locsol}).
The proof also indicates that the Parodi's relation plays an essential role in obtaining proper higher-order energy inequalities for the E--L system.
Next, we address a Beale--Kato--Majda type criterion that characterizes the first
finite singular time of the local classical solution in terms of $\|\curlu\|_{L^\infty}$ as well as $\|\nabla \bd\|_{L^\infty}$ (cf. Theorem \ref{blow}).
For the uncoupled incompressible Navier--Stokes equations (or Euler equation), in their pioneering work \cite{BKM}  Beale, Kato and Majda showed the criterion in
terms of $L^\infty$-norm of the vorticity $\curlu$.
As far as liquid crystal systems are concerned, the situation is more involved due to interactions between the fluids and molecules.
In the recent work \cite{WZZ12} the authors obtain a BKM type criterion for the full E--L system with constraint $|\bd|=1$ in terms of $\|\curlu\|_{L^\infty}$
and $\|\nabla \bd\|_{L^\infty}$ (see also \cite{HW12} for a simplified system).
For our system \eqref{e1}--\eqref{e3} the main difficulty comes from the lack of control on the term $\|\bd\|_{L^\infty}$, which brings troubles
to estimate the nonlinear kinematic transport terms as well as the higher-order Leslie stress tensors.
\beos
With minor modifications of the proof of Section 5, we can obtain existence of global weak solutions also for the case of the homogeneous Neumann boundary
condition for the director field, i.e., $\partial_{\mathbf{n}} \bd = \mathbf{0}$, which is suitable for the implementation of a numerical scheme (cf. \cite{LS}).
The case of periodic boundary conditions for $\ub$ and $\bd$ can be treated with even simpler computations.
However, when we consider classical solutions (and their blow-up criterion), we are only able to handle periodic boundary
conditions. The main technical reason is that the periodic boundary conditions allow us to integrate by parts to obtain proper  higher-order
differential energy inequalities (cf. also \cite{sunliu, wuxuliu}.)
\eddos

\paragraph{Plan of the paper.} The remaining part of this paper is organized as follows. Section 2 contains the notations and the constraints on the
Leslie coefficients. In Section 3, we first introduce a proper weak formulation of the E--L system and state the weak solutions existence theorem. After obtaining some useful (lower-order) \textit{a priori} estimates from the basic energy law, we prove the existence of weak solutions with finite energy.
Finally, Section 4 is devoted to the proof of local well-posedness of local classical solutions and a BKM type criterion.


\section{Preliminaries and assumptions}
\label{sec:mainres}

Let $\Omega\subset \RR^3$ be a bounded domain with  smooth boundary $\Gamma$.
We denote by $L^p(\Omega;\RR)$ and $W^{m,p}(\Omega;\RR)$ the usual $L^p$-spaces and Sobolev spaces of real measurable functions on $\Omega$ and by
$L^p(\Omega;\RR^3)$, $L^p(\Omega;\RR^{3\times 3})$, $W^{m,p}(\Omega;\RR^3)$, $W^{m,p}(\Omega;\RR^{3\times3})$ the corresponding spaces of vector functions.
Sometimes they will be simply denoted by $L^p, W^{m,p}$.
If $I$ is an interval of $\RR^+$ and $X$ a Banach space, we also use the function space $L^p(I;X)$, $1 \leq p \leq +\infty$, which consists of $p$-integrable
functions with values in $X$. Next, we introduce some function spaces related to the Dirichlet boundary value problem  (cf. \cite{LinLiusimply,temam}):
 \begin{eqnarray}
 && \mathcal{V}=C_0^\infty(\Omega;\mathbb{R}^3)\cap \{\textbf{v}: \nabla\cdot \textbf{v}=0\},\no\\
 && \textbf{H}=\ \text{the closure of}\ \mathcal{V} \text{\ in\ } L^2(\Omega;
\mathbb{R}^3),\no\\
 && \textbf{V}=\ \text{the closure of\ } \mathcal{V} \text{\ in\ } \ W^{1,2}_0(\Omega; \RR^3),\no\\
 && \textbf{V}'=\text{the\ dual of\ } \textbf{V}.\no
 \end{eqnarray}
When we consider the periodic boundary value problem, we set $\Omega=Q$, where $Q$ is a unit box in $\RR^3$.
Then we denote by $H^m_p(Q)$, $m\in \mathbb{N}$, the space of functions which are in $W^{m,2}_{loc}(\mathbb{R}^3)$ and periodic in space with period $Q$.
For the sake of simplicity, we denote the inner product on $L^2$-spaces by $(\cdot,\cdot)$ and the associated norm by $\|\cdot\|$.
Throughout this paper, the same letter $C$ stands for a constant which may be different each time it appears.

We consider suitable conditions on the physical coefficients $\lambda_1, \lambda_2, \mu_1, ..., \mu_6$,
so that the E--L system \eqref{e1}--\eqref{e3}, together with proper boundary conditions, obeys certain dissipative
energy law (cf. \cite{LL01}). More precisely, we take the constraints
 \begin{eqnarray}
 && \lambda_1<0, \label{lama1a}\\
 && \mu_1\geq 0, \quad \mu_4>0, \label{mu14}\\
 && \mu_5+\mu_6 \geq 0,  \label{mu56}  \\
 &&\lambda_1=\mu_2-\mu_3, \ \ \ \lambda_2=\mu_5-\mu_6, \label{lama1}\\
&&\mu_2+\mu_3=\mu_6-\mu_5. \label{lam2}
\end{eqnarray}
Assumptions \eqref{lama1a}--\eqref{mu14} are made to provide necessary conditions for the dissipation of the system
(cf. \cite{E91, Le68, Le79}).
Relations \eqref{lama1} are necessary conditions to satisfy the equation of
motion identically (cf. \cite[Section 6]{Le68}).


Equation \eqref{lam2}
is called \emph{Parodi's relation}, which is derived from Onsager reciprocal relations expressing the equality
between flows and forces in thermodynamic systems out of equilibrium  (cf. \cite{P70}).
The Parodi's relation yields a constraint on the Leslie coefficients such that the dynamics of an incompressible nematic
liquid crystal flow now involves five independent coefficients in \eqref{v5} instead of six. In \cite{WXL12},
the authors show that the Parodi's relation also serves as a stability condition for the E--L system \eqref{e1}--\eqref{e3}.

In this paper, we will assume the following two different sets of hypotheses on the coefficients, as proposed in \cite{WXL12},
\begin{itemize}
\item \textbf{Case 1} (with Parodi's relation). Suppose that
 \eqref{lama1a}--\eqref{lam2} are satisfied. Moreover, we assume
 \begin{equation} \frac{(\lambda_2)^2}{-\lambda_1} \leq
 \mu_5+\mu_6.
 \label{critical point of lambda 2}
 \end{equation}
 \item \textbf{Case 2} (without Parodi's relation). Suppose that \eqref{lama1a}--\eqref{lama1} are satisfied. Moreover, we assume
 \begin{equation}
|\lambda_2-\mu_2-\mu_3| <
2\sqrt{-\lambda_1}\sqrt{\mu_5+\mu_6}.
\label{noPa1}
  \end{equation}

\end{itemize}

As far as the potential $W$ is concerned, for the sake of simplicity, we will always assume that it satisfies the
 Ginzburg--Landau approximation, that is
\begin{equation}
W(\bd)=\frac{1}{4\varepsilon^2}(|\bd|^2-1)^2.\label{GL1}
\end{equation}

\beos
In general, $W$ may be written as a sum of a convex part and a smooth, but possibly non-convex one, for instance,
\begin{align}
&W \in C^2(\RR^3), \quad W \geq 0,\no\\
&W=W_1+W_2 \, \hbox{ s.t. } \,  W_1 \hbox{ is convex and } W_2 \in C^1(\RR^3), \, \nabla W_2 \in C^{0,1}(\RR^3;\mathbb{R}^3)\no.
\end{align}
With minor modifications in the proof as in \cite{CR}, it is easy to check that our result on existence of weak solutions to the
E--L system \eqref{e1}--\eqref{e3} (cf. Theorem \ref{theo1} below) also holds for the general potential $W$ given above.
However, in order to prove existence of (local) classical solution (cf. Theorem \ref{locsol}), additional assumptions on smoothness
of the potential $W$ should be added.
\eddos

\section{Existence of weak solutions}

First we state the weak formulation of the problem under consideration, in which the momentum equation (\ref{e1}) together
with the incompressibility condition (\ref{e2}) are replaced by integral identities.
 \bedef \label{defweak} {\rm (Weak solution).}
A pair $(\ub,\,\bd)$ with $\ub(0,\cdot)=\ub_0$,  $\bd(0, \cdot)=\bd_0$, a.e. in $\Omega$, is called a {\em weak solution}
to the E--L system \eqref{e1}--\eqref{e3} subject to the Dirichlet boundary conditions \eqref{slip} and \eqref{diri}, if it
 belongs to the class
\begin{equation} \label{reg1}
\ub \in L^\infty(0,T; L^2(\Omega; \RR^3))
\cap L^2(0,T;W^{1,2}_0(\Omega; \RR^3)),
\end{equation}
\begin{equation}\label{reg1bis}
\partial_t\ub \in L^2(0,T;W^{-1,\frac65}(\Omega; \RR^3)),
\end{equation}
\begin{equation} \label{reg2}
\bd \in L^\infty(0,T;
W^{1,2}(\Omega; \RR^3))\cap L^2(0,T; W^{2,2}(\Omega;\RR^3))\cap H^1(0,T; L^{\frac32}(\Omega;\RR^3)),
\end{equation}
and, moreover,
\begin{equation}
\label{weak1}
\int_{\Omega} \ub(t, \cdot) \cdot \Grad \varphi = 0, \quad \mbox{for a.e.} \ t \in (0,T),
\end{equation}
for any  $\varphi\in W^{1,2}(\Omega;\RR)$,
\begin{eqnarray}
 && \langle \partial_t\ub, \mathbf{v} \rangle_{W^{-1, \frac65}, W_0^{1,6}}
- \int_{\Omega}\ub\otimes\ub :\Grad\mathbf{v} dx \nonumber\\
 &=& \int_{\Omega} (\Grad \bd\odot\Grad \bd) : \Grad \mathbf{v}dx- \io \bsig : \nabla \mathbf{v}dx,\quad \mbox{for a.e.} \ t \in (0,T), \label{weak2}
\end{eqnarray}
for any $\mathbf{v} \in W^{1,6}_0({\Omega}; \RR^3)$ such that  $\nabla \cdot \mathbf{v}=0$.
Finally, the equation for the molecule director is satisfied in the strong sense, that is
\begin{eqnarray}
\label{weak3}
\bd_t+(\ub\cdot\nabla )\bd-\omega
\bd+\frac{\lambda_2}{\lambda_1}A\bd&=&-\frac{1}{\lambda_1}\left(\Delta
\bd-\Grad_{\bd} W(\bd)\right), \,
\mbox{ a.e. in }(0,T) \times \Omega,\\
\label{weak4}
\bd &=& \bd_0|_\Gamma,
\hskip3truecm \ \mbox{a.e. on }(0,T)\times\Gamma.
\end{eqnarray}
\eddef

\beos We note that the choice of the exponent $\frac65$ in \eqref{reg1bis} is due to the regularity of the nonlinear term of highest order in the stress tensor
(see \eqref{ap7a} below). Besides, a weak solution $(\ub,\,\bd)$ corresponding to the Neumann or periodic boundary conditions can be defined in a similar way.
\eddos

Introduce now the total energy of the system \eqref{e1}--\eqref{e3}, consisting of kinetic and potential energies, given by
\begin{equation} \mathcal{E}=\frac{1}{2}\|\ub\|^2+\frac{1}{2}\|\nabla \bd\|^2+\int_{\Omega}W(\bd)dx.
\label{total energy of the system}
\end{equation}

The main result of this section read as follows and its proof will be postponed to Subsection 3.2.

\bete \label{theo1} {\rm (Existence of weak solution with finite energy).} Let $\Omega \subset \RR^3$ be a bounded domain of class $C^{2}$.
Assume (\ref{GL1}) and the initial data such that
\begin{eqnarray}
&&\ub_0 \in L^{2}(\Omega; \RR^3),\quad  \dive \ub_0 = 0 \ \mbox{in} \,  L^2(\Omega),\label{hyp4}\color{black} \\
&&\bd_0 \in W^{1,2} (\Omega; \RR^3), \quad \bd_0|_\Gamma\in W^{\frac32,2}(\Gamma; \RR^3).\label{hyp5}
\end{eqnarray}
Then, for both \textbf{Case 1} and \textbf{Case 2}, the Dirichlet problem (\ref{weak1})--(\ref{weak4}) possesses a global-in-time weak
solution $(\ub, \bd)$, in the sense of Definition \ref{defweak}.
Moreover, for a.e.~$t\in(0,T),$ the following energy inequalities hold, namely, for \textbf{Case 1},
\begin{eqnarray}
&& \mathcal{E}(t)+ \int_0^t \left[\io \Big(\mu_1|\bd^{T}A\bd|^2+\frac{\mu_4}{2}|\nabla \ub|^2
\Big)dx -\frac{1}{\lambda_1}\|\Delta \bd-\Grad_{\bd} W(\bd)\|^2\right]d\tau
 \nonumber\\
 &&\ \ +\Big[\mu_5+\mu_6+\frac{(\lambda_2)^2}{\lambda_1}\Big]\int_0^t\|A\bd\|^2 d\tau \leq \mathcal{E}(0),\no
 \end{eqnarray}
while for \textbf{Case 2}
\begin{eqnarray}
\mathcal{E}(t) +\int_0^t \left[\io \left(\mu_1|\bd^{T}A\bd|^2+\frac{\mu_4}{2}|\nabla \ub|^2
\right)dx +\eta (\|A\bd\|^2+\|\mathcal{N}\|^2)\right] d\tau \leq \mathcal{E}(0).\no
\end{eqnarray}
\ente

\beos
Note that the assumption $\mathbf{d}_0\in W^{\frac32, 2}(\Gamma)$  is required to derive the regularity $\bd\in L^2(0,T; W^{2,2}(\Omega;\RR^3))$ stated in \eqref{reg2} by means of the classical elliptic estimate (cf.~also \eqref{apr33}) .
\eddos


\subsection{\textit{A priori} estimates}

We establish here a number of formal a priori estimates. These will assume a rigorous character in the framework of the
approximation scheme presented in Subsection 3.2 below.

Let $(\ub, \bd)$ be a solution to the E--L system  \eqref{e1}--\eqref{e3} subject to one of the following boundary conditions:
(i) $\ub|_\Gamma=\mathbf{0}$, $\bd|_{\Gamma}=\bd_0$; (ii) $\ub|_\Gamma=\mathbf{0}$, $\partial_\mathbf{n}\bd|_{\Gamma}=\mathbf{0}$;
or (iii) $\Omega=Q$ with $\ub$, $\bd$ periodic in space.
Following the detailed calculations as in \cite[Lemma 2.1, Lemma 2.2]{WXL12} and \cite[Theorem 2.1]{LL01}, for
all the three types of boundary conditions for $\bd$  (i.e., Dirichlet, Neumann and periodic), we are able to obtain
some formal energy inequality for the E--L system \eqref{e1}--\eqref{e3}.
More precisely, we have
\begin{eqnarray}
\frac{d}{dt}\mathcal{E}(t)&=& -\int_{\Omega}(\mu_1|\bd^{T}A\bd|^2+\frac{\mu_4}{2}|\nabla \ub|^2+(\mu_5+\mu_6)|A\bd|^2)dx \nonumber\\
&&+\lambda_1\|\mathcal{N}\|^2+(\lambda_2-\mu_2-\mu_3)(\mathcal{N}, A\bd).
\label{BELnoPa}
\end{eqnarray}
As a consequence,
\begin{itemize}
\item  \textbf{Case 1}. The total energy $\mathcal{E}(t)$ is decreasing in time and it holds
 \begin{eqnarray}
 \frac{d}{dt}\mathcal{E}(t)
&=& -\io \Big(\mu_1|\bd^{T}A\bd|^2+\frac{\mu_4}{2}|\nabla \ub|^2
\Big)dx +\frac{1}{\lambda_1}\|\Delta \bd-\Grad_{\bd} W(\bd)\|^2
 \nonumber\\
 &&-\Big(\mu_5+\mu_6+\frac{(\lambda_2)^2}{\lambda_1}\Big)\|A\bd\|^2\leq 0.
 \label{basic energy law at the critical point of lambda 2}
 \end{eqnarray}

\item  \textbf{Case 2}. The total energy $\mathcal{E}(t)$  is decreasing in time and there exists a small constant
  $\eta>0$ such that
  \begin{eqnarray}
 \frac{d}{dt}\mathcal{E}(t)  &=&-\int_{\Omega}(\mu_1|\bd^{T}A\bd|^2+\frac{\mu_4}{2}|\nabla
\ub|^2+(\mu_5+\mu_6)|A\bd|^2)dx \nonumber\\
&&+\lambda_1\|\mathcal{N}\|^2+(\lambda_2-\mu_2-\mu_3)(\mathcal{N}, A\bd)\nonumber\\
&\leq&
-\io \left(\mu_1|\bd^{T}A\bd|^2+\frac{\mu_4}{2}|\nabla \ub|^2
\right)dx -\eta (\|A\bd\|^2+\|\mathcal{N}\|^2)\nonumber\\
&\leq& 0.\label{belc}
 \end{eqnarray}
 \end{itemize}

 \beos\label{eqqq}
 We can easily see from the equation \eqref{e3} that the energy dissipations in \textbf{Case 1} and \textbf{Case 2}
 are indeed equivalent. In particular, recalling also \eqref{def1}, we have
 $$ \|A\bd\|^2+\|\Delta \bd-\Grad_{\bd} W(\bd)\|^2\approx \|A\bd\|^2+\|\mathcal{N}\|^2.$$
  \eddos

From now on we just present the case in which $(\ub, \bd)$ satisfies the boundary conditions of Dirichlet type.
Corresponding results for Neumann or periodic boundary conditions can be obtained in a similar way.

From the previous energy laws we can obtain now (formal) \textit{a priori} energy estimates for the solutions to the E--L system
\begin{eqnarray}
&&\|\ub\|_{L^\infty(0,T; L^2(\Omega;\RR^3))}\leq C,\quad   \|\ub\|_{L^2(0,T;W^{1,2}_0(\Omega; \RR^3))}\leq C,
\label{eapr1}\\
&&\|\bd\|_{L^\infty(0,T; W^{1,2}(\Omega;\RR^3))}\leq C,
\label{eapr2}\\
&&
\|-\Delta \bd +\Grad_{\bd} W(\bd)\|_{L^2(0,T; L^2(\Omega;\RR^3))}\leq C,\label{eapr2bis}\\
&& \|A\bd\|_{L^2(0,T; L^2(\Omega;\RR^3))}\leq C,\quad \|\mathcal{N}\|_{L^2(0,T; L^2(\Omega;\RR^3))}\leq C,\label{ANes}\\
&& \|\bd^{T}A\bd\|_{L^2(0,T; L^2(\Omega;\RR^3))}\leq C,\label{dAdes}
\end{eqnarray}
where $C$ is a positive constant depending on $\|\ub_0\|$, $\|\bd_0\|_{W^{1,2}}$, $\Omega$, and the
coefficients of the system, but it is independent of $T$.

Hence, estimates \eqref{eapr1}--\eqref{dAdes} imply that
\begin{eqnarray}
\ub &\in& L^\infty(0,T; L^2(\Omega;\RR^3))\cap L^2(0,T;W^{1,2}_0(\Omega; \RR^3)),
\label{apr1}\\
\bd &\in& L^\infty(0,T; W^{1,2}(\Omega;\RR^3)),
\label{apr2}\\
-\Delta \bd +\Grad_{\bd} W(\bd) &\in& L^2((0,T)\times \Omega; \RR^3),\label{apr2bis}\\
A\bd, \, \mathcal{N} &\in& L^2(0,T; L^2(\Omega;\RR^3)),
\label{apr22}\\
 \bd^{T}A\bd &\in& L^2(0,T; L^2(\Omega;\RR)).
\end{eqnarray}
It follows from \eqref{apr1} that  $\ub\otimes \ub\in L^2(0,T;L^{\frac32}(\Omega; \RR^{3\times 3}))$.
 Then, using the incompressible condition \eqref{e2}, we get
  \begin{equation}
  (\ub\cdot\nabla) \ub=\nabla \cdot (\ub\otimes \ub)\in L^2(0,T;W^{-1,\frac32}(\Omega; \RR^3)).
  \end{equation}
From the definition of $W$ and \eqref{apr2},
we have
\begin{equation}\label{apr3}
\nabla_{\bd} W(\bd)\in L^2(0,T;L^2(\Omega; \RR^3)).
\end{equation}
 By the classical elliptic regularity theory, property \eqref{apr2bis} together with  \eqref{apr3} and the Dirichlet boundary condition \eqref{diri} for $\bd$ yield
\begin{equation} \label{apr33}
\bd\in L^2 (0,T; W^{2,2}(\Omega; \RR^3)),
\end{equation}
where the boundedness constant also depends on $\|\bd_0\|_{W^{\frac32,2}(\Gamma;\mathbb{R}^3)}$.
As a consequence, we have
$$(\ub\cdot\Grad)\bd, \ \omega\bd\  \in \ L^2(0,T;L^{\frac32}(\Omega; \RR^3)).$$
Then, from \eqref{apr22} and the definition of $\mathcal{N}$, we deduce
\begin{equation}\label{apr4}
\partial_t \bd\in L^2(0,T;L^{\frac32}(\Omega;\RR^3)).
\end{equation}
We recall the well-known anisotropic Sobolev embedding theorem for the Banach space
$$ V_2((0,T)\times \Omega)=L^\infty(0, T; L^2(\Omega))\cap L^2(0, T; W^{1,2}(\Omega))$$
such that, when the spatial dimension is three, it holds
$ V_2((0,T)\times \Omega)\hookrightarrow L^{\frac{10}{3}}((0,T)\times \Omega)$.
Therefore,  from \eqref{apr2} and \eqref{apr33} we infer
\begin{equation}\label{apr5}
\Grad\vc{d} \in L^{\frac{10}{3}}(0,T;L^{\frac{10}{3}}(\Omega;\RR^{3\times3})).
\end{equation}
In a similar manner, we get
\begin{equation}\label{apr5u}
\ub \in L^{\frac{10}{3}}(0,T;L^{\frac{10}{3}}(\Omega;\RR^{3})).
\end{equation}
Besides, using the fact that, for  $s\in [0,1]$, it holds (see e.g., \cite[Definition 1.1, pp. 27]{LP} for the definition of the interpolation space $(\cdot, \cdot)_s$)
$$
\big(L^2(0,T; W^{2,2}(\Omega; \mathbb{R}^3)), L^\infty(0,T; W^{1,2}(\Omega; \mathbb{R}^3))\big)_s= L^\frac{2}{1-s}(0,T; W^{2-s,2}(\Omega; \mathbb{R}^3)),
$$
and the Sobolev embedding in $3D$ such that $W^{2-s,2}(\Omega; \mathbb{R}^3)\hookrightarrow L^\frac{6}{2s-1}(\Omega; \mathbb{R}^3)$,
then, taking $s=\frac45$, we have (cf. \cite{RS04})
$$
\big(L^2(0,T; W^{2,2}(\Omega; \mathbb{R}^3)), L^\infty(0,T; W^{1,2}(\Omega; \mathbb{R}^3))\big)_s\hookrightarrow L^{10}(0,T; L^{10}(\Omega; \mathbb{R}^3)).
$$
As a result,  from \eqref{apr2} and \eqref{apr33} we infer
\begin{equation}
\bd\in L^{10}(0,T; L^{10}(\Omega; \mathbb{R}^3)).\label{d10}
\end{equation}
From the above estimates \eqref{ANes}, \eqref{dAdes}, \eqref{apr5}, \eqref{d10} and proper interpolation inequalities, we see that the nonlinear
stress terms fulfil
 \begin{eqnarray}
  \nabla \bd\odot \nabla \bd, \,
    A\bd \otimes \bd, \, \mathcal{N}\otimes \bd
   &\in& L^2(0,T;L^{\frac32}(\Omega; \RR^{3\times 3}))\cap L^{\frac53}((0,T)\times\Omega;\RR^{3\times3}),
   \label{ap7}\\
    (\bd^{T}A\bd)\bd\otimes\bd &\in& L^2(0,T;L^{\frac65}(\Omega; \RR^{3\times 3}))\cap L^{\frac{10}{7}}((0,T)\times\Omega;\RR^{3\times3}). \label{ap7a}
 \end{eqnarray}
 For instance, for the nonlinear stress tensor of the highest order (cf.  \eqref{ap7a}), we have used the following facts
 \begin{eqnarray}
  \|(\bd^{T}A\bd)\bd\otimes\bd\|_{L^2(0,T;L^{\frac65}(\Omega; \RR^{3\times 3}))}&\leq& \|\bd^{T}A\bd\|_{L^2(0,T;L^{2}(\Omega; \RR))}\|\bd\|^2_{L^\infty(0,T;L^{6}(\Omega; \RR^{3}))}\no\\
  &\leq& C\|\bd^{T}A\bd\|_{L^2(0,T;L^{2}(\Omega; \RR))}\|\bd\|^2_{L^\infty(0,T;W^{1,2}(\Omega; \RR^{3}))}\no
 \end{eqnarray}
 and (cf. \eqref{d10})
 \begin{eqnarray}
 \|(\bd^{T}A\bd)\bd\otimes\bd\|_{L^{\frac{10}{7}}((0,T)\times\Omega;\RR^{3\times3})}&\leq& \|\bd^{T}A\bd\|_{L^2(0,T;L^{2}(\Omega; \RR))}\|\bd\|_{L^{10}(0,T; L^{10}(\Omega; \mathbb{R}^3))}^2.\no
 \end{eqnarray}
Finally, since
  $$L^2(0,T;L^{\frac32}(\Omega; \RR^{3\times 3}))\subset L^2(0,T;L^{\frac65}(\Omega; \RR^{3\times 3})),$$
then the (distributional) divergence of the Leslie stress tensor $\bsig$ satisfies
 $$\grad \cdot \bsig \in L^2(0,T;W^{-1,\frac65}(\Omega; \RR^3)).$$

\subsection{Proof of Theorem \ref{theo1}}
 The proof of Theorem \ref{theo1} consists of several steps. As in \cite{CR,ffrs}, we shall construct a suitable family of
\emph{approximate} problems whose solutions weakly converge (always up to subsequences)
to certain limit functions that solve the problem in the sense of Definition \ref{defweak}.
We note that in order to prove the existence of weak solutions, the double-approximation scheme described below is not necessary for the special and simpler case when the molecule director $\mathbf{d}$ obeys a weak maximum principle, see, for instance, \cite{LinLiusimply, LL01}.

\textbf{Step 1. Approximation.}
We introduce a double-approximation scheme that consists of a standard Faedo--Galerkin method with an
approximation of the convective term as well as a regularization of the momentum equation (\ref{weak1})--(\ref{weak2}) by adding an
$r$-Laplacian operator acting on the velocity $\ub$.

We take the orthonormal basis $\{ \mathbf{v}_n \}_{n=1}^\infty$ of the Hilbert space
$\mathbf{V}$. Fixing $M, \,N \in \mathbb{N}$ such that $M\leq N$, we consider the finite-dimensional space
$X_N = {\rm span}\{ \mathbf{v}_n\}_{n=1}^N$.
 The approximate velocity field $\ub_{N,M} \in C^1([0,T]; X_N)$ solves
the following equations
$$
 \frac{{\rm d}}{{\rm d}t} \int_\Omega \ub_{N,M} \cdot \mathbf{v} dx-
\int_{\Omega} [ \ub_{N,M} ]_M\otimes  \ub_{N,M} : \Grad \mathbf{v}dx + \frac{1}{M}\io |\Grad\ub_{N,M}|^{r-2}\Grad\ub_{N,M}\cdot\Grad\vb dx
$$
\begin{equation}
 = \int_\Omega  \Grad \bd_{N,M} \odot \Grad \bd_{N,M} : \Grad \mathbf{v}dx- \io \bsig_{N,M} : \nabla \mathbf{v}dx,
 \quad \hbox{ for all }t \in (0,T), \label{approx1}
\end{equation}
\begin{equation}
\label{approx1bis}
\int_\Omega \ub_{N,M}(0, \cdot) \cdot \mathbf{v} dx = \int_\Omega
\ub_0 \cdot \mathbf{v}dx,
\end{equation}
for any $\mathbf{v} \in X_N$ and certain fixed $r\in (\frac{10}{3}, +\infty)$. Here the symbol $[ \mathbf{v} ]_M$ denotes
the orthogonal projection onto the finite-dimensional space $X_M = {\rm span}\{\mathbf{v}_n\}_{n=1}^M$.
The $r$-Laplacian regularization term $\frac{1}{M}|\Grad\ub_{N,M}|^{r-2}\Grad\ub_{N,M}$ in \eqref{approx1}
is introduced to obtain enough regularity for the velocity field. This enables us to deduce certain energy inequalities for the
limit functions $(\ub_M, \bd_M)$, after we pass to the limit $N\to +\infty$ (see Step 2 below).
The approximate function $\bd_{N,M}$ for the molecule director is determined in terms of $\ub_{N,M}$ as the unique solution of the
parabolic system
\begin{equation}
\partial_t \bd_{N,M} + \ub_{N,M} \cdot \Grad \bd_{N,M} -\omega_{N,M}
\bd_{N,M}+\frac{\lambda_2}{\lambda_1}A_{N,M}\bd_{N,M}\no
\end{equation}
\begin{equation}
= \Delta \bd_{N,M}-
\nabla_{\bd} W(\bd_{N,M}), \quad\quad \hbox{ in $(0,T) \times \Omega$},\label{approx2}
\end{equation}
\begin{eqnarray}
\bd_{N,M} &=& \bd_0|_\Gamma,\quad \hbox{ on $(0,T) \times \Gamma$},\label{approx3}
\\
\bd_{N,M}(0, \cdot) &=& \bd_{0}, \quad \hbox{ in $\Omega$},\label{approx4}
\end{eqnarray}
where
$$
 A_{N,M}=\frac12(\nabla \ub_{N,M}+\nabla^{T}\ub_{N,M}),\
\,\; \omega_{N,M}=\frac12(\nabla \ub_{N,M}-\nabla^{T}\ub_{N,M}).$$

For any fixed $M$ and $N$, we can solve problem (\ref{approx1})--(\ref{approx4}) by means of a fixed-point argument as in \cite[Chapter 3]{FN}
(cf. also \cite{LinLiusimply} for a simplified model of the E--L system).
Indeed, we observe that all the \emph{a priori} bounds derived formally in the previous section still hold for the approximate problem.
Hence, if we fix $\tilde{\ub} \in C([0,T]; X_N)$, then we can find
$\bd = \bd[\tilde{\ub}]$ solving (\ref{approx2})--(\ref{approx4}). Inserting $\bd[\tilde{\ub}]$ in system
(\ref{approx1})--(\ref{approx1bis}),  we can define a mapping $\tilde{\ub} \mapsto {\cal T}[\tilde{\ub}]$, where $\ub={\cal T}[\tilde{\ub}]$
is the solution to the system.
On account of the \emph{a priori} bounds, we can easily show that ${\cal T}$ admits a fixed point by means of the classical Schauder's argument
on $(0,T_0)$, with $0< T_0 \leq T$.
Finally, since the \emph{a priori} estimates are independent of the final time $T_0$, we are allowed to conclude that the approximate solutions can be
extended to the whole time interval $[0,T]$. The details are omitted here.\medskip

\textbf{Step 2: Passage to the limit as $N\to +\infty$.} Hereafter we just treat \textbf{Case 1} since \textbf{Case 2} can be handled in a similar way (cf. Remark \ref{eqqq}).
From \eqref{approx1}--\eqref{approx1bis}, we deduce that the regular approximate solutions $(\ub_{N,M}, \bd_{N,M})$ satisfy
\begin{eqnarray}\label{est}
&& \frac{1}{2}\|\ub_{N,M}(t)\|^2+\frac{1}{2}\|\nabla
\bd_{N,M}(t)\|^2+\int_{\Omega}W(\bd_{N,M}(t))dx\no\\
&&\ \  + \int_0^t \io \Big(\mu_1|\bd^{T}A_{N,M}\bd_{N,M}|^2+\frac{\mu_4}{2}|\nabla \ub_{N,M}|^2
+\frac{1}{M} |\Grad\ub_{N,M}|^r \Big)dxd\tau
 \nonumber\\
 &&\ \ +\frac{1}{\lambda_1}\int_0^t \|\Delta \bd_{N,M}-\Grad_{\bd} W(\bd_{N,M})\|^2 d\tau\no\\
 &&\ \ +\Big[\mu_5+\mu_6+\frac{(\lambda_2)^2}{\lambda_1}\Big]\int_0^t\|A_{N,M}\bd_{N,M}\|^2 d\tau \no\\
 &=& \frac{1}{2}\|\ub_{N,M}(0)\|^2+\frac{1}{2}\|\nabla
\bd_{N,M}(0)\|^2+\int_{\Omega}W(\bd_{N,M}(0))dx\label{appener}
 \end{eqnarray}
 and, from the formal estimates in the previous section,
 that the following convergence results hold (up to a subsequence)
\begin{align}
\label{cuN}
&\ub_{N,M}\to\ub_M \ \mbox{ weakly-(*) in } L^\infty(0,T;L^2(\Omega;\RR^3))
\cap L^2(0,T;W^{1,2}(\Omega;\RR^3)),\\
\label{cdN}
&\bd_{N,M}\to \bd_M \ \mbox{ weakly-(*) in } L^\infty(0,T;W^{1,2}(\Omega;\RR^3))
\cap L^2(0,T;W^{2,2}(\Omega;\RR^3)),\\
&\dt \bd_{N,M}\to \dt \bd_M \ \mbox{ weakly in } L^2(0,T;L^{\frac32}(\Omega;\RR^3)).\label{cdtNbis}
\end{align}
By virtue of \eqref{cdN}, \eqref{cdtNbis} and the Aubin--Lions lemma (cf. e.g., \cite{simon}), we
have, for any arbitrary small $\xi$ such that $0< \xi <<1$,
\begin{equation}
\bd_{N,M}\to \bd_M \ \mbox{ strongly
in } C([0,T]; W^{1-\xi, 2}(\Omega; \RR^{3}))\cap L^2(0, T; W^{2-\xi, 2}(\Omega; \RR^{3})).\label{scond}
\end{equation}
 Then a simple interpolation argument yields that
\begin{equation}\label{custgd}
\bd_{N,M}\to \bd_M \ \mbox{ strongly in }L^\eta((0,T)\times\Omega; \RR^{3}),\ \mbox{ for }\eta\in \left(1,10\right),
\end{equation}
\begin{equation}\label{custg}
\Grad\bd_{N,M}\to \Grad \bd_M \ \mbox{ strongly
in }L^\eta((0,T)\times\Omega; \RR^{3\times3}),\ \mbox{ for }\eta_1\in \left(1,\frac{10}{3}\right).
\end{equation}
On account of the regularizing $r$-Laplacian term introduced in \eqref{approx1}, we obtain the following additional regularity
from the corresponding energy estimate
\begin{equation}\label{unuova}
  \frac{1}{M}\|\Grad\ub_{N,M}\|^r_{L^r((0,T)\times\Omega;\mathbb{R}^{3\times 3})}\leq C,
  \ \mbox{ for } r \in \left(\frac{10}{3}, +\infty\right),
\end{equation}
where the positive constant $C$ does not depends on $M$ and $N$. Hence,
 \begin{equation}\label{nonnu}
 \frac{1}{M}|\Grad\ub_{N,M}|^{r-2}\Grad\ub_{N,M}\in L^{\frac{r}{r-1}}((0,T)\times\Omega;\mathbb{R}^{3\times3}),
 \ \ \text{with}\ \frac{r}{r-1}\in \left(1, \frac{10}{7}\right).
 \end{equation}
Estimate \eqref{unuova} implies that
\begin{align}
\label{cuNnuova}
&\Grad\ub_{N,M}\to\Grad\ub_M \ \mbox{ weakly in } L^{r}(0,T;L^r(\Omega;\RR^{3\times 3})).
\end{align}
Furthermore, we have (recall also \eqref{ap7}--\eqref{ap7a})
\begin{equation}
\label{cutN}
 \dt\ub_{N,M}\to\dt\ub_M \ \mbox{ weakly in }
  L^{\frac{r}{r-1}}(0,T; W^{-1, \frac{r}{r-1}}(\Omega;\RR^3))+L^2(0,T; W^{-1,2}(\Omega;\RR^3)).
\end{equation}
  Combining  \eqref{cuNnuova} with \eqref{cutN} and using the Aubin--Lions lemma once more, we get
  \begin{equation}\label{ustro}
  \ub_{N,M}\to\ub_M \ \mbox{ strongly in } L^r(0,T;W^{1-\xi, r}(\Omega;\RR^3)), \quad \forall \, \xi \, \, \,\text{s.t.} \,\,  \, 0<\xi<<1,
  \end{equation}
  which together with the Sobolev embedding yields (recall that $r>\frac{10}{3}>3$)
  \begin{equation}
  \ub_{N,M}\to\ub_M \ \mbox{ strongly in } L^r(0,T;L^{\infty}(\Omega;\RR^3)). \label{ustroa}
  \end{equation}
Then interpolating  \eqref{cuN} with \eqref{ustroa}, we obtain
\begin{equation}\label{custrN}
  \ub_{N,M}\to\ub_M \ \mbox{ strongly in } L^{r+2}((0,T)\times\Omega;\RR^3).
\end{equation}
A combination of \eqref{custrN}  with \eqref{custg} gives
\begin{equation}\label{convls}
 \ub_{N,M}\cdot\Grad\bd_{N,M}\to \ub_M\cdot\Grad \bd_M\
 \mbox{ strongly in } L^{s_1}((0,T)\times\Omega;\RR^3), \ \mbox{with}\   \frac{1}{s_1}=\frac{1}{r+2}+\frac{1}{\eta_1}.
\end{equation}
Let $r=\frac{10}{3}+\varsigma$, $\varsigma>0$ and take $\eta_1=\frac{10}{3}-\varsigma'$ for certain small $0<\varsigma'\leq\min\{\varsigma,1\}$.
It is easy to see that
\begin{equation}
s_1= \frac{\left(\frac{16}{3}+\varsigma\right)\left(\frac{10}{3}-\varsigma'\right)}{\frac{26}{3}+\varsigma-\varsigma'}> \frac{80}{39} -\left(\frac{8}{13}+\frac{3}{26}\varsigma\right)\varsigma'>2,
\end{equation}
provided that $0<\varsigma'<<1$.
Similarly, combining \eqref{cuNnuova} and \eqref{custgd}, one obtains, for some $ s_2>2$,
\begin{equation}
A_{N,M} \bd_{N,M} \to A_M \bd_M , \quad \omega_{N,M} \bd_{N,M}\to \omega_M \bd_M,\  \mbox{ weakly in } L^{s_2}((0,T)\times \Omega;\RR^3),\label{cAd}
\end{equation}
where
$$
 A_{M}=\frac12(\nabla \ub_{M}+\nabla^{T}\ub_{M}),\
\,\; \omega_{M}=\frac12(\nabla \ub_{M}-\nabla^{T}\ub_{M}).$$
The above relations \eqref{convls}, \eqref{cAd} and the $L^2$-bound of $-\Delta \bd_{N,M}+ \nabla_{\bd}W(\bd_{N,M})$
(cf. estimate \eqref{est}) imply that
\begin{equation}
\dt \bd_{N,M}\to \dt \bd_M \ \mbox{ weakly in } L^2(0,T;L^{2}(\Omega;\RR^3)).\label{cdtNbisa}
\end{equation}
Finally, from the strong convergence \eqref{scond} we easily see that (cf. assumptions on $W$)
\begin{equation}
\nabla_{\bd}W(\bd_{N,M})\to \nabla_{\bd}W(\bd_{M})\quad \mbox{strongly in } L^2((0,T)\times \Omega;\RR^3).\label{conW}
\end{equation}
As a consequence, we are able to pass to the limit as $N \to \infty$ in the approximate equation \eqref{approx2} to get,
a.e. in $(0,T)\times\Omega$,
\begin{equation} \label{aapprox3}
\partial_t \bd_{M}
+ \ub_{M} \cdot \Grad \bd_{M} -\omega_{M}
\bd_{M}+\frac{\lambda_2}{\lambda_1}A_{M}\bd_{M}
 = \Delta \bd_{M}-\nabla_{\bd}W(\bd_{M}),
\end{equation}
with
\begin{eqnarray}
\bd_{M} &=& \mathbf{d}_0, \quad \hbox{a.e.~in } (0,T) \times \Gamma, \label{aapprox4}\\
\bd_{M}(0, \cdot) &=& \bd_{0}|_\Gamma, \quad \hbox{a.e.~in } \Omega.\label{aapprox5}
\end{eqnarray}

In order to pass to the limit $N\to+\infty$ in the approximate momentum equation \eqref{approx1}, we investigate the nonlinear terms.
First, it follows from \eqref{custrN} that
\begin{equation}\label{custrNa}
  [\ub_{N,M}]_M\otimes \ub_{N,M} \to[\ub_M]_M\otimes\ub_M \ \mbox{ strongly in } L^{\frac{r+2}{2}}((0,T)\times\Omega;\RR^3),
\end{equation}
while \eqref{custg} yields, for $\eta_2\in \left(1,\frac{5}{3}\right)$,
\begin{equation}\label{custg1}
\Grad\bd_{N,M} \odot \Grad\bd_{N,M}  \to \Grad \bd_M \odot \Grad \bd_M\ \mbox{ strongly
in } L^{\eta_2}((0,T)\times\Omega; \RR^{3\times3}).
\end{equation}
Next, from \eqref{scond} and \eqref{cAd}, we have
\begin{eqnarray}
&& A_{N,M}\bd_{N,M}\otimes \bd_{N,M} \to A_{M}\bd_{M}\otimes \bd_{M}\ \mbox{ weakly in } L^2(0,T; L^\frac32(\Omega; \RR^{3\times 3})),\no\\
&&  \omega_{N,M}\bd_{N,M}\otimes \bd_{N,M} \to \omega_{M}\bd_{M}\otimes \bd_{M}\ \mbox{ weakly in } L^2(0,T; L^\frac32(\Omega; \RR^{3\times 3})),\no
\end{eqnarray}
while it follows from \eqref{scond}, \eqref{convls} and \eqref{cdtNbisa} that
\begin{eqnarray}
&& (\ub_{N,M}\cdot \nabla \bd_{N,M})\otimes \bd_{N,M} \to (\ub_{M}\cdot \nabla \bd_{M})\otimes \bd_{M}\
\mbox{ strongly in } L^2(0,T; L^\frac32(\Omega; \RR^{3\times 3})),\no\\
&& \dt \bd_{N,M}\otimes \bd_{N,M} \to \dt \bd_{M}\otimes \bd_{M}\ \mbox{ weakly in } L^2(0,T; L^{\eta_3}(\Omega; \RR^{3\times 3})),\
\eta_3\in \left(1, \frac32\right).\no
\end{eqnarray}
Concerning the stress term of the highest-order $(\bd^T_{N,M}A_{N,M} \bd_{N,M})\bd_{N,M}\otimes \bd_{N,M}$,
from \eqref{cuN} and \eqref{custgd} we first infer that, for some $ s_3\in (1, 10/7)$,
\begin{equation}
\bd_{N,M}^TA_{N,M}\bd_{N,M} \to \bd_M^T A_M\bd_M \ \mbox{ weakly in } L^{s_3}((0,T)\times \Omega;\RR).
\end{equation}
This can be improved due to the $L^2$-estimate of $\bd_{N,M}^TA_{N,M}\bd_{N,M}$ (cf. \eqref{appener}) and uniqueness of the weak limit
\begin{equation}
\bd_{N,M}^TA_{N,M}\bd_{N,M} \to \bd_M^T A_M\bd_M \ \mbox{ weakly in } L^{2}((0,T)\times \Omega;\RR).
\end{equation}
Combining it with \eqref{custgd} again, we have, for some $s_4\in (1, 10/7)$,
\begin{eqnarray}
&& (\bd^T_{N,M}A_{N,M} \bd_{N,M})\bd_{N,M}\otimes \bd_{N,M}\no \to  (\bd^T_MA_M \bd_M)\bd_M\otimes \bd_M,\\
&& \mbox{ weakly in } L^{s_4}((0,T)\times \Omega; \RR^{3\times 3}). \label{highest}
\end{eqnarray}
Besides,  estimate \eqref{nonnu} implies
\begin{equation}
|\Grad\ub_{N,M}|^{r-2}\Grad \ub_{N,M}\to \overline{|\Grad\ub_{M}|^{r-2}\Grad \ub_{M}},\ \
\mbox{weakly in } L^{\frac{r}{r-1}}((0,T)\times \Omega; \RR^{3\times3}).\label{mity}
\end{equation}
Then the pair of limit functions $(\ub_M, \bd_M)$ solves the problem (note that in \eqref{approx1} the projection on $X_M$
is kept in the convective term, cf. \eqref{custrNa})
\begin{equation}
\int_\Omega \ub_M\cdot \nabla \varphi dx=0, \quad \mbox{for a.e.} \ t\in (0,T),\label{aapprox0}
\end{equation}
for any $\varphi\in C^\infty(\bar \Omega)$, and
$$
 \int_0^t \langle\partial_t \ub_{M},  \mathbf{v} \rangle d\tau-\int_0^t
\int_{\Omega} [ \ub_{M} ]_M\otimes  \ub_{M} : \Grad \mathbf{v}dx d\tau+ \frac{1}{M}\int_0^t\io \overline{|\Grad\ub_{M}|^{r-2}\Grad \ub_{M}}:\Grad\vb dxd\tau
$$
\begin{equation}
 = \int_0^t\int_\Omega  \Grad \bd_{M} \odot \Grad \bd_{M} : \Grad \mathbf{v}dxd\tau-\int_0^t \io \bsig_{M} : \nabla \mathbf{v}dxd\tau,
 \quad \hbox{ for all }t \in (0,T), \label{aapprox1}
\end{equation}
for any $\mathbf{v} \in C^\infty(\overline{\Omega}; \RR^3)$ with $\nabla \cdot \mathbf{v}  = 0$. Here
\begin{eqnarray}
\bsig_M&=&\mu_1(\bd^T_MA_M \bd_M)\bd_M\otimes \bd_M+\mu_2\mathcal{N}_M\otimes \bd_M+\mu_3\bd_M\otimes \mathcal{N}_M+
\mu_4A_M\no\\
&& +\mu_5(A_M\bd_M)\otimes \bd_M+\mu_6\bd_M\otimes(A_M\bd_M)\no\\
&\in& L^\frac{10}{7}(0,T;L^\frac{10}{7}(\Omega;\RR^{3\times3})). \label{v5M}
\end{eqnarray}
It remains to show that in \eqref{aapprox1}
\begin{equation}
|\Grad\ub_{M}|^{r-2}\Grad \ub_{M}=\overline{|\Grad\ub_{M}|^{r-2}\Grad \ub_{M}}.\label{kkk}
\end{equation}
On account of the $L^r$-regularity ($r>\frac{10}{3}$) of $\Grad\ub_M$ (cf. \eqref{cuNnuova}) and \eqref{v5M},
we are allowed to take $\mathbf{v}=\ub_M$ as a test function in \eqref{aapprox1}. On the other hand, we multiply \eqref{aapprox3}
by $\Delta\bd_{M} -\nabla_{\bd} W(\bd_{M})$ and integrate on $(0, t) \times \Omega$ (this is possible since
\eqref{aapprox3} makes sense in $L^2(Q_T)$, cf.~(\ref{convls})--(\ref{cdtNbisa})). Adding the two resultants together, we get
\begin{eqnarray}
&& \frac{1}{2}\|\ub_{M}(t)\|^2+\frac{1}{2}\|\nabla
\bd_{M}(t)\|^2+\int_{\Omega}W(\bd_{M}(t))dx\no\\
&& \ \ + \int_0^t \io \left(\mu_1|\bd^{T}_MA_{M}\bd_{M}|^2+\frac{\mu_4}{2}|\nabla \ub_{M}|^2\right)dxd\tau
\no\\
&& \ \ +\frac{1}{M} \int_0^t \io \overline{|\Grad\ub_{M}|^{r-2}\Grad \ub_{M}} :\Grad\ub_{M} dxd\tau
 \nonumber\\
 && \ \ +\frac{1}{\lambda_1}\int_0^t\io \|\Delta \bd_{M}-\Grad_{\bd} W(\bd_{M})\|^2 dxd\tau\no\\
 && \ \ +\Big(\mu_5+\mu_6+\frac{(\lambda_2)^2}{\lambda_1}\Big)\int_0^t\|A_{M}\bd_{M}\|^2 d\tau \no\\
 &=& \frac{1}{2}\|\ub_{M}(0)\|^2+\frac{1}{2}\|\nabla
\bd_{M}(0)\|^2+\int_{\Omega}W(\bd_{M}(0))dx.
\label{appener1}
 \end{eqnarray}
Passing to the limit as $N\to+\infty$ in \eqref{appener} and using the lower semi-continuity of norms, we infer from \eqref{appener1} that
\[
\limsup_{N\to+\infty} \int_0^T\io |\Grad \ub_{N,M}|^r dxdt\leq \int_0^T\io \overline{|\Grad\ub_{M}|^{r-2}\Grad\ub_M}:\Grad\ub_{M}dxdt
\]
Recalling \eqref{cuNnuova} and \eqref{mity}, by means of the well-known Minty's trick for monotone operators (cf. \cite{M2} or \cite[Lemma 3.2.2]{Z04}), we deduce that \eqref{kkk} holds.
This concludes the passage to the limit as $N\to +\infty$.\medskip


\textbf{Step 3: Passage to the limit as $M\to+\infty$.} The final step is to pass to the limit as
$M\to+\infty$ in (\ref{aapprox3})--(\ref{aapprox5}) and \eqref{aapprox0}--\eqref{v5M}.
First, we observe that we are able to obtain the same convergence results as in (\ref{cuN})--\eqref{custg}, \eqref{conW}, \eqref{custg1}
and \eqref{highest} for $\ub_M$ and $\bd_M$, while letting
$M\to +\infty$, since the $r$-Laplacian regularization is not necessary to perform the limit in these terms.
Moreover, we deduce from \eqref{kkk} and  \eqref{appener1} that
\begin{equation}
  \frac{1}{M}\|\Grad\ub_{M}\|^r_{L^r((0,T)\times\Omega;\mathbb{R}^{3\times 3})}\leq C, \quad r \in \left(\frac{10}{3}, +\infty\right),\label{nuMe}
\end{equation}
where $C$ is independent of $M$. As a consequence, the following convergence result holds true (keeping in mind that $\frac{r}{r-1}<\frac{10}{7}$)
\begin{align}
\label{cdtuM}
& \dt\ub_M\to\dt\ub\ \mbox{ weakly in } L^{\frac{r}{r-1}}(0,T;W^{-1,\frac{r}{r-1}}(\Omega;\RR^3)).
\end{align}
Using \eqref{cuN} and the Aubin--Lions lemma again, we get
\begin{equation}
\ub_{M}\to\ub \ \mbox{ strongly in } L^2(0,T;W^{1-\xi,2}(\Omega;\RR^3)).
\end{equation}
Then we conclude that,  for some  $ s_5\in(1, 5/4)$ and  $s_6\in (1, 5/3)$,
\begin{equation}\label{convlsM}
 \ub_{M}\cdot\Grad\bd_{M}\to \ub \cdot\Grad \bd\
 \mbox{ strongly in } L^{s_5}((0,T)\times\Omega;\RR^3),
\end{equation}
\begin{equation}
A_{M} \bd_{M} \to A \bd , \quad \omega_{M} \bd_{M}\to \omega \bd,\  \mbox{ weakly in } L^{s_6}((0,T)\times \Omega;\RR^3).\label{cAdM}
\end{equation}
Next, from \eqref{custgd} (for $\bd_M$ instead of $\bd_{N,M}$) and \eqref{cAdM}, we deduce
\begin{eqnarray}
&& A_{M}\bd_{M}\otimes \bd_{M} \to A\bd\otimes \bd\ \mbox{ weakly in } L^{s_7}((0,T)\times\Omega; \RR^{3\times 3}),\no\\
&&   \omega_{M}\bd_{M}\otimes \bd_{M} \to \omega\bd\otimes \bd\ \mbox{ weakly in } L^{s_7}((0,T)\times\Omega; \RR^{3\times 3}),\no
\end{eqnarray}
 with $s_7\in (1, 10/7)$.
Besides, for some $s_8\in (1, 10/9)$ and $s_{9}\in (1, 6/5)$, we have
\begin{eqnarray}
&& (\ub_{M}\cdot \nabla \bd_{M})\otimes \bd_{M} \to (\ub\cdot \nabla \bd)\otimes \bd\ \mbox{ strongly in } L^{s_8}((0,T) \times \Omega; \RR^{3\times 3}),\no\\
&& \dt \bd_{M}\otimes \bd_{M} \to \dt \bd\otimes \bd\ \mbox{ weakly in } L^2(0,T; L^{s_{9}}(\Omega; \RR^{3\times 3})).\no
\end{eqnarray}
\color{black}
Finally, concerning the $r$-Laplacian regularization term, by the interpolation inequality and \eqref{nuMe}, we obtain
 \begin{eqnarray}
 && \left\|M^{-\frac{1}{r-1}}\Grad\ub_{M}\right\|
 _{L^{r-1}((0,T)\times\Omega;\mathbb{R}^{3\times3})}
 \no\\
 &\leq& M^{-\frac{1}{(r-1)(r-2)}}\left( M^{-\frac{1}{r}}\|\nabla \ub_M\|_{L^{r}((0,T)\times\Omega;\mathbb{R}^{3\times3})}\right)^a
 \|\nabla \ub_M\|_{L^{2}((0,T)\times\Omega;\mathbb{R}^{3\times3})}^{1-a}
  \no\\
  &\leq& CM^{-\frac{1}{(r-1)(r-2)}},\quad \mbox{with}\ \  a=\frac{r(r-3)}{(r-1)(r-2)}\in (0, 1),\label{limM}
 \end{eqnarray}
 which yields
\begin{equation}
\label{cduMM}
M^{-\frac{1}{r-1}}\Grad \ub_{M}\to 0 \mbox{ strongly in } L^{r - 1}((0,T)\times \Omega;\RR^{3\times3}).
\end{equation}
We are now in a position to pass to the limit as $M \to + \infty$ in (\ref{aapprox3})--(\ref{aapprox5}) and \eqref{aapprox0}--\eqref{v5M},
and finally recover the system (\ref{weak1})--(\ref{weak4}). The proof of Theorem \ref{theo1} is complete.


\section{Local well-posedness and BKM type blow-up criterion}
\label{BLUP}

In this section we first prove the existence and uniqueness of local classical solutions to the E--L system \eqref{e1}--\eqref{e3} subject to
periodic boundary conditions and then establish a BKM type blow-up criterion.
For this purpose, the following results will be frequently used in the subsequent proofs

\bele \label{KP} {\rm (cf. \cite{KP}).} For $s\geq 0$, there holds
\begin{eqnarray}
&& \|fg\|_{H^s}\leq C(\|f\|_{L^\infty}\|g\|_{H^s}+\|f\|_{H^s}\|g\|_{L^\infty}),\no\\
&& \|\Lambda^s (fg)-f\Lambda^s g\|\leq C(\|f\|_{H^s}\|g\|_{L^\infty}+\|\nabla f\|_{L^\infty}\|g\|_{H^{s-1}}),\no
\end{eqnarray}
 where $\Lambda$ is the Fourier multiplier such that  $\Lambda f(x)=\sum_{k\in
 \mathbb{Z}^d}(1+|k|^2)^\frac{s}{2}e^{2\pi ik\cdot x}\widehat{f}(k)$.
\enle
\bele \label{tr} {\rm (cf. \cite{Tr}).} For any $s>0$, we denote by $[s]$ the integer part of $s$.
Assume that $F(\cdot)$ is a smooth function on $\mathbb{R}$ with $F(0)=0$
and $f \in H^s(\Omega)\cap L^\infty (\Omega)$.
Then we have
$$ \|F(f)\|_{H^s}\leq c(1+\|f\|_{L^\infty})^{[ s] +1}\|f\|_{H^s}, $$
where the constant $c$ depends on $\displaystyle\sup_{0\leq k\leq [ s]+2, \ |y|\leq \|f\|_{L^\infty}}\|F^{(k)}(y)\|_{L^\infty}$.
\enle

\subsection{Local classical solutions}
\bete
 \label{locsol} {\rm (Local classical solution).} Let $\Omega:=Q=(0,1)^3$ be the unit cubic in $\mathbb{R}^3$.
 Suppose that the conditions in \textbf{Case 1} are satisfied and $s\geq 3$ is an integer.
 For any $\ub_0\in H^s_p(Q)$ with $\nabla \cdot \ub_0=0$ and $\bd_0\in H^{s+1}_p(Q)$, there exists a $T_0>0$ depending on $\|\ub_0\|_{H^3}$, $\|\bd_0\|_{H^4}$
 and the coefficients of the system, but uniform in $s$, such that the E--L system \eqref{e1}--\eqref{e3} subject to periodic boundary conditions admits a unique
local classical solution satisfying
\begin{align}
& \ub \in C([0, T_0]; H^s_p(Q))\cap L^2(0,
T_0; H^{s+1}_p(Q)),\label{s1a}\\
&\bd \in C([0, T_0]; H^{s+1}_p(Q))\cap
L^2(0, T_0; H^{s+2}_p(Q)).\label{s2a}
\end{align}

\ente

\begin{proof}
In order to prove the existence and uniqueness of local classical solutions, we can first construct approximate solutions to the periodic value problem of the
E--L system \eqref{e1}--\eqref{e3} by the Galerkin method (e.g., \cite{WZ12}). We note that here we no longer need the double-approximation scheme described in Section 3, since we can make use of higher-order estimates. After obtaining suitable uniform estimates on the approximate solutions,
we are able to pass to the limit. Since the approximation procedure is standard, here below we just perform the necessary uniform  higher-order energy estimates.

Under the assumptions of \textbf{Case 1}, we can see that the (lower-order) uniform estimates \eqref{eapr1}--\eqref{dAdes} are still satisfied. Moreover, recalling \cite[Lemma 5.1]{WXL12},
the quantity $$\mathcal{A}(t)=\|\nabla \ub\|^2+\|\Delta \bd- \nabla_\bd W(\bd)\|^2$$ satisfies the following differential inequality:
\begin{equation}
\frac{d}{dt}\mathcal{A}(t)\leq C(
\mathcal{A}(t)^6+\mathcal{A}(t)),\label{A}
\end{equation}
where $C$ is a constant that only depends on the coefficients $\mu_i$,
$\lambda_i$, $\|\ub_0\|$ and $\|\bd_0\|_{H^1}$. It follows from \eqref{A} that there exists a certain $T_1<+\infty$ such that
\begin{equation}
\mathcal{A}(t)\leq C, \quad t\in [0,T_1],
\end{equation}
for some constant $C$ depending only on  $\|\ub_0\|_{H^1}$, $\|\bd_0\|_{H^2}$ and the coefficients of the system. As a result
\begin{equation}
\sup_{0\leq t\leq T_1} \|\ub(t)\|_{H^1}+\|\bd(t)\|_{H^2}\leq C.\label{pes1}
\end{equation}
In particular, by the Sobolev embedding theorem, we also have
 \begin{equation}
 \|\bd(t)\|_{L^\infty}\leq C, \quad \forall\ t\in [0,T_1]. \label{linfinity}
 \end{equation}

On account of estimates \eqref{pes1} and \eqref{linfinity}, we proceed to obtain higher-order estimates for $\ub$ and $ \bd$.
For $s\geq 3$, applying $\Lambda^s$ to \eqref{e1} and testing the resultant by $\Lambda^s\ub $, we get
\begin{eqnarray}
&& \frac12\frac{d}{dt}\|\Lambda^s \ub\|^2+ \frac{\mu_4}{2}\|\nabla \Lambda^s \ub\|^2\no\\
&=& -(\Lambda^s(\ub\cdot\nabla \ub), \Lambda^s \ub)+ (\Lambda^s(\nabla \bd \odot \nabla \bd ), \nabla \Lambda^s \ub)\no\\
&& -\mu_1(\Lambda^s[(\bd^TA \bd)\bd\otimes \bd], \nabla \Lambda^s \ub)\no\\
&& -\mu_2(\Lambda^s(\mathcal{N}\otimes \bd), \nabla\Lambda^s \ub)
-\mu_3(\Lambda^s(\bd\otimes \mathcal{N}), \nabla\Lambda^s \ub)\no\\
&& -\mu_5(\Lambda^s[(A\bd)\otimes \bd], \nabla\Lambda^s \ub)
-\mu_6(\Lambda^s[\bd\otimes(A\bd)], \nabla\Lambda^s \ub)\no\\
&:=& \sum_{i=1}^7 K_i.\no
\end{eqnarray}
 Moreover, applying $\nabla \Lambda^s$ to \eqref{e3} and testing the resultant by $\nabla \Lambda^s\bd$, we obtain
 \begin{eqnarray}
 && \frac12\frac{d}{dt}\|\nabla \Lambda^s \bd\|^2 -\frac{1}{\lambda_1}\|\Delta \Lambda^s \bd\|^2\no\\
 &=&
 -(\nabla \Lambda^s(\ub\cdot\nabla \bd), \nabla \Lambda^s \bd)-(\Lambda^s(\omega
\bd), \Delta \Lambda^s \bd) \no\\
&&+\frac{\lambda_2}{\lambda_1}(\Lambda^s(A\bd), \Delta \Lambda^s \bd)+\frac{1}{\lambda_1}(\Lambda^s[\Grad_{\bd} W(\bd)], \Delta \Lambda^s \bd)\no\\
&:=& \sum_{i=8}^{11} K_i.\no
 \end{eqnarray}
For $K_1$, using the incompressibility of $\ub$, we have
\begin{eqnarray}
K_1&=& -(\Lambda^s(\ub\cdot\nabla \ub)-\ub\cdot \nabla \Lambda^s \ub, \Lambda^s \ub)\no\\
&\leq& C\|\Lambda^s(\ub\cdot\nabla \ub)-\ub\cdot \Lambda^s \nabla \ub\|\|\Lambda^s \ub\|\no\\
&\leq& C\|\nabla \ub\|_{L^\infty}\|\ub\|_{H^s}^2.\no
\end{eqnarray}
The term $K_2$ can be simply estimated as
\begin{eqnarray}
K_2&\leq& \|\Lambda^s(\nabla \bd\odot\nabla \bd)\|\|\nabla \Lambda^s \ub\|\no\\
 &\leq& \kappa \|\nabla \Lambda^s \ub\|^2+C\|\nabla \bd\|_{L^\infty}^2\|\nabla \Lambda^s \bd\|^2.\no
\end{eqnarray}
Similarly to $K_1$, for $K_8$, we get
\begin{eqnarray}
K_8&=&  -(\nabla \Lambda^s(\ub\cdot\nabla \bd)-\ub\cdot\nabla (\nabla \Lambda^s \bd), \nabla \Lambda^s \bd)\no\\
&\leq &  C(\|\ub\|_{H^{s+1}}\|\nabla \bd\|_{L^\infty}+\|\nabla \ub\|_{L^\infty}\|\nabla \bd\|_{H^{s}})\|\nabla \Lambda^s \bd\|\no\\
&\leq& \kappa \|\nabla \Lambda^s \ub\|^2+ C(1+\|\nabla \ub\|_{L^\infty}+\|\nabla \bd\|_{L^\infty}^2)\|\nabla \Lambda^s \bd\|^2.\no
\end{eqnarray}
By the interpolation inequality
$$ \|\nabla \ub\|_{L^\infty}\leq C\|\ub\|_{H^1}^\frac{2s-5}{2s-2}\|\Lambda^s\ub\|^\frac{3}{2s-2},\quad s\geq 3,$$
and since $\frac{3}{2s-2}\leq \frac12$ when $s\geq 4$, we have
\begin{eqnarray}
&& \|\bd\|_{H^s}^2\|A\|_{L^\infty}^2\leq C\|\bd\|_{H^2}\|\bd\|_{H^4}\|\nabla \ub \|_{L^\infty}^\frac23\|\ub\|_{H^1}^\frac13\|\Lambda^3\ub\|,\quad s=3, \no\\
&& \|\bd\|_{H^s}^2\|A\|_{L^\infty}^2\leq C\|\bd\|_{H^{s-1}}\|\bd\|_{H^{s+1}}\|\ub\|_{H^1}\|\Lambda^s\ub\|, \quad s\geq 4.\no
\end{eqnarray}
From \eqref{pes1} we infer, for $s\geq 3$,
\begin{equation}
\|\bd\|_{H^s}^2\|A\|_{L^\infty}^2\leq C(1+\|\nabla \ub\|^\frac23_{L^\infty})(\|\Lambda^s\ub\|^2+ \|\bd\|_{H^{s-1}}^2\|\bd\|_{H^{s+1}}^2), \quad t\in [0,T_1].\label{EE}
\end{equation}
Then, for the stress term $K_3$, noticing that $\omega$ is antisymmetric and using \eqref{pes1}, \eqref{EE}, we obtain
\begin{eqnarray}
K_3
&=& -\mu_1(\Lambda^s[(\bd^TA \bd)\bd\otimes \bd],  \Lambda^s(A+\omega))\no\\
&=& -\mu_1(\Lambda^s[(\bd^TA \bd)\bd\otimes \bd],  \Lambda^3 A)\no\\
&=& -\mu_1 \|\bd^T \Lambda^s A \bd\|^2+ \mu_1(\Lambda^s[(\bd^TA \bd)\bd\otimes \bd]- [(\bd^T \Lambda^s A \bd)\bd\otimes \bd], \Lambda^s A)
\no\\
&\leq& -\mu_1 \|\bd^T  \Lambda^s A \bd\|^2 +\kappa \|\nabla \Lambda^s \ub\|^2 + C\|\Lambda^s[(\bd^TA \bd)\bd\otimes \bd]- [(\bd^T \Lambda^s A \bd)\bd\otimes \bd]\|^2 \no\\
&\leq& -\mu_1 \|\bd^T  \Lambda^s A \bd\|^2 +\kappa \|\nabla \Lambda^s \ub\|^2 \no\\
&&
+C(\|\bd\|_{L^\infty}^3\|\bd\|_{H^s}\|A\|_{L^\infty}+\|\bd\|_{L^\infty}^3\|\nabla \bd\|_{L^\infty}\|A\|_{H^{s-1}})^2
\no\\
&\leq& -\mu_1 \|\bd^T \Lambda^s A \bd\|^2 +\kappa \|\nabla \Lambda^s \ub\|^2 +C(1+\|\nabla \ub\|_{L^\infty})(\|\Lambda^s\ub\|^2+ \|\bd\|_{H^{s-1}}^2\|\bd\|_{H^{s+1}}^2)\no\\
&& + C\|\nabla \bd\|_{L^\infty}^2\|\Lambda^s \ub\|^2
\no\\
&\leq& -\mu_1 \|\bd^T \Lambda^s A \bd\|^2 +\kappa \|\nabla \Lambda^s \ub\|^2\no\\
&& + C(1+\|\nabla \ub\|_{L^\infty}+\|\nabla \bd\|_{L^\infty}^2)(\|\Lambda^s \ub\|^2+\|\Lambda^{s-1}\bd\|^2\|\Lambda^{s+1} \bd\|^2).\no
\end{eqnarray}
Using the definition of $A, \omega$, and the fact that $\omega$ is antisymmetric, then it holds
\begin{eqnarray}
&& K_4+K_5+K_6+K_7\no\\
&=& -(\mu_2+\mu_3)(\Lambda^s(\mathcal{N}\otimes \bd), \Lambda^s A)-(\mu_2-\mu_3)(\Lambda^s(\mathcal{N}\otimes \bd), \Lambda^s \omega)\no\\
&&-(\mu_5+\mu_6)(\Lambda^s((A\bd)\otimes \bd), \Lambda^s A)-(\mu_5-\mu_6)(\Lambda^s((A\bd)\otimes \bd), \Lambda^s \omega)\no\\
&:=& P_1+P_2+P_3+P_4.\no
\end{eqnarray}
On the other hand, using the equation  \eqref{e3} such that $\lambda_1\mathcal{N}+\lambda_2 A\bd=-\Delta \bd+\nabla_\bd W(\bd)$, we have
\begin{eqnarray}
&& K_9+K_{10}\no\\
&=&
  \lambda_1(\Lambda^s \mathcal{N}, \Lambda^s(\omega
\bd))
  +\lambda_2(\Lambda^s (A\bd), \Lambda^s(\omega
\bd))\no\\
&& -\lambda_2(\Lambda^s \mathcal{N} , \Lambda^s (A\bd))
  -\frac{(\lambda_2)^2}{\lambda_1}(\Lambda^s (A\bd), \Lambda^s (A\bd))\no\\
&& -(\Lambda^s (\Grad_\bd W(\bd)), \Lambda^s(\omega
\bd))
   + \frac{\lambda_2}{\lambda_1} (\Lambda^s (\Grad_\bd W(\bd)), \Lambda^s (A\bd))\no\\
&:=& \sum_{i=5}^{10} P_i.\no
\end{eqnarray}
Recalling the condition \eqref{lama1} and the Parodi's relation \eqref{lam2}, we can recombine the above terms $P_1,...,P_{10}$ and eliminate the
term with highest-order derivative as follows
\begin{eqnarray}
&& P_1+P_7+P_{10}\no\\
&=& \lambda_2 (\Lambda^s(\mathcal{N}\otimes \bd), \Lambda^s A)-\lambda_2(\Lambda^s \mathcal{N} , \Lambda^s (A\bd))+\frac{\lambda_2}{\lambda_1} (\Lambda^s (\Grad_\bd W(\bd)), \Lambda^s (A\bd))\no\\
&=& \lambda_2 (\Lambda^s(\mathcal{N}\otimes \bd)-\Lambda^s \mathcal{N}\otimes \bd, \Lambda^s A)- \lambda_2(\Lambda^s \mathcal{N} , \Lambda^s (A\bd)-\Lambda^s A \bd)\no\\
&& +\frac{\lambda_2}{\lambda_1} (\Lambda^s (\Grad_\bd W(\bd)), \Lambda^s (A\bd))\no\\
&=& \lambda_2 (\Lambda^s(\mathcal{N}\otimes \bd)-\Lambda^s \mathcal{N}\otimes \bd, \Lambda^s A)
+\frac{(\lambda_2)^2}{\lambda_1} (\Lambda^s (A\bd), \Lambda^s (A\bd)-\Lambda^s A \bd)\no\\
&& +\frac{\lambda_2}{\lambda_1}(\Lambda^s \Delta \bd, \Lambda^s (A\bd)-\Lambda^s A \bd)+\frac{\lambda_2}{\lambda_1}(\Lambda^s \Grad_\bd W(\bd), \Lambda^s A \bd)\no\\
&:=& R_1+R_2+R_3+R_4.\no
\end{eqnarray}
\begin{eqnarray}
&& P_2+P_5+P_9\no\\
&=& \lambda_1(\Lambda^s \mathcal{N}, \Lambda^s(\omega
\bd))-\lambda_1(\Lambda^s(\mathcal{N}\otimes \bd), \Lambda^s \omega)\no\\
&& -(\Lambda^s (\Grad_\bd W(\bd)), \Lambda^s(\omega
\bd))\no\\
&=& \lambda_1(\Lambda^s \mathcal{N}, \Lambda^s(\omega
\bd)-\Lambda^s\omega \bd)-\lambda_1(\Lambda^s(\mathcal{N}\otimes \bd)-\Lambda^s \mathcal{N}\otimes \bd, \Lambda^s \omega)\no\\
&& -(\Lambda^s (\Grad_\bd W(\bd)), \Lambda^s(\omega
\bd))\no\\
&=& -\lambda_1(\Lambda^s(\mathcal{N}\otimes \bd)-\Lambda^s \mathcal{N}\otimes \bd, \Lambda^s \omega)
-\lambda_2(\Lambda^s (A\bd), \Lambda^s(\omega
\bd)-\Lambda^s\omega \bd)\no\\
&& -(\Lambda^s \Delta \bd, \Lambda^s(\omega
\bd)-\Lambda^s\omega \bd)+(\Lambda^s (\Grad_\bd W(\bd)),\Lambda^s\omega \bd)\no\\
&:=& R_5+R_6+R_7+R_8.\no
\end{eqnarray}
\begin{eqnarray}
P_4+P_6&=& -\lambda_2(\Lambda^s((A\bd)\otimes \bd), \Lambda^s \omega)+\lambda_2(\Lambda^s (A\bd), \Lambda^s(\omega
\bd))\no\\
&=& -\lambda_2(\Lambda^s((A\bd)\otimes \bd)-\Lambda^s (A \bd) \otimes \bd, \Lambda^s \omega)\no\\
&& +\lambda_2(\Lambda^s (A\bd), \Lambda^s(\omega
\bd)-\Lambda^s \omega \bd)\no\\
&:=& R_9+R_{10}.\no
\end{eqnarray}
\begin{eqnarray}
P_3+P_8&=& -(\mu_5+\mu_6)(\Lambda^s((A\bd)\otimes \bd), \Lambda^s A) -\frac{(\lambda_2)^2}{\lambda_1}(\Lambda^s (A\bd), \Lambda^s (A\bd))\no\\
&=& -\left(\mu_5+\mu_6+\frac{(\lambda_2)^2}{\lambda_1}\right)(\Lambda^s (A\bd), \Lambda^s (A\bd))\no\\
&&  -(\mu_5+\mu_6)(\Lambda^s((A\bd)\otimes \bd)-\Lambda^s(A\bd)\otimes \bd, \Lambda^s A)\no\\
&&  +(\mu_5+\mu_6)(\Lambda^s(A\bd), \Lambda^s( A \bd) -\Lambda^s A \bd)\no\\
&:=& R_{11}+R_{12}+R_{13}.\no
\end{eqnarray}
As a consequence, in order to estimate the terms $K_4, K_5, K_6, K_7, K_9, K_{10}$, we can turn to estimate the new terms $R_1,...,R_{13}$.

By the following bounds
\begin{eqnarray}
&& \|\bd\|_{H^3}\leq C\|\bd\|_{H^{4}}^\frac13\|\nabla \bd\|_{L^\infty}^\frac23, \no\\
&& \|\Delta \bd\|_{L^\infty}\leq C\|\bd\|_{H^3}^\frac12\|\bd\|_{H^4}^\frac12,\no
\end{eqnarray}
   we deduce, for $ t\in [0,T_1]$,
\begin{eqnarray}
\|\bd\|_{H^s}^2\|\Delta \bd\|_{L^\infty}^2
&\leq& C\|\bd\|_{H^3}^3\|\bd\|_{H^4}\leq C\|\nabla \bd\|_{L^\infty}^2\|\bd\|_{H^4}^2,\ \ \text{when}\ s=3, \label{EE2}
\\
\|\bd\|_{H^s}^2\|\Delta \bd\|_{L^\infty}^2&\leq& C\|\bd\|_{H^{s-1}}\|\bd\|_{H^{s+1}}\|\bd\|_{H^3}\|\bd\|_{H^4}\no\\
&\leq& C\|\bd\|_{H^{s-1}}^2\|\bd\|_{H^{s+1}}^2, \quad \text{when} \ s\geq 4.\label{EE3}
\end{eqnarray}
Next, using equation \eqref{e3}, Lemmas \ref{KP}, \ref{tr} and estimates \eqref{pes1}--\eqref{EE3}, we have
\begin{eqnarray}
R_1 &\leq& \kappa\|\Lambda^s A\|^2+ C(\|\bd\|_{H^s}\|\mathcal{N}\|_{L^\infty}+\|\nabla \bd\|_{L^\infty}\|\mathcal{N}\|_{H^{s-1}})^2\no\\
&\leq & \kappa\|\Lambda^3 A\|^2+C \|\bd\|_{H^s}^2(\|A\|_{L^\infty}^2\|\bd\|_{L^\infty}^2+\|\Delta \bd\|_{L^\infty}^2+\|\Grad_\bd W(\bd)\|_{L^\infty}^2)\no\\
&& +C \|\nabla \bd\|_{L^\infty}^2(\|A\|_{H^{s-1}}^2\|\bd\|_{L^\infty}^2+\|A\|_{L^\infty}^2\|\bd\|_{H^{s-1}}^2+\|\Delta \bd\|^2_{H^{s-1}}+\|\Grad_\bd W(\bd)\|_{H^{s-1}}^2)\no\\
&\leq& \kappa\|\Lambda^3 A\|^2 +C(1+\|\nabla \ub\|_{L^\infty}+ \|\nabla \bd\|_{L^\infty}^2)(\|\Lambda^s\ub\|^2+ \|\bd\|_{H^{s-1}}^2\|\bd\|_{H^{s+1}}^2)\no\\
&&
 + C\|\bd\|_{H^{s-1}}\|\bd\|_{H^{s+1}}(\|\Delta \bd\|_{H^1}\|\Delta \bd\|_{H^2}+ 1)\no\\
&&+C \|\nabla \bd\|_{L^\infty}^2(\|\Lambda^s \ub\|^2+\|\Lambda^{s-1}\bd\|^2+\|\Lambda^{s-1}\bd\|^4+\|\Lambda^{s+1}\bd\|^2)\no\\
&\leq& \kappa\|\Lambda^3 A\|^2+ C(1+\|\nabla \ub\|_{L^\infty}+\|\nabla \bd\|_{L^\infty}^2)[\|\Lambda^s\ub\|^2+
(1+\|\bd\|_{H^{s-1}}^2)\|\bd\|_{H^{s+1}}^2].\no
\end{eqnarray}
\begin{eqnarray}
R_2&\leq& C(\|A\|_{H^s}\|\bd\|_{L^\infty}+\|A\|_{L^\infty}\|\bd\|_{H^s})(\|\bd\|_{H^s}\|A\|_{L^\infty}+\|\nabla \bd\|_{L^\infty}\|A\|_{H^{s-1}})\no\\
&\leq& \kappa\|\Lambda^s A\|^2+C \|\bd\|_{H^s}^2\|\nabla \ub\|_{L^\infty}^2+C\|\nabla \bd\|_{L^\infty}^2\|\Lambda^s \ub\|^2\no\\
&\leq& \kappa\|\Lambda^s A\|^2+ C(1+\|\nabla \ub\|_{L^\infty}+\|\nabla \bd\|_{L^\infty}^2)(\|\nabla \Lambda^s \bd\|^2+\|\Lambda^s\ub\|^2).\no
\end{eqnarray}
\begin{eqnarray}
R_3&\leq& \kappa\|\Delta \Lambda^s \bd\|^2+C(\|\bd\|_{H^s}\|A\|_{L^\infty}+\|\nabla \bd\|_{L^\infty}\|A\|_{H^{s-1}})^2\no\\
&\leq&  \kappa\|\Lambda^s A\|^2+ C(1+\|\nabla \ub\|_{L^\infty}+\|\nabla \bd\|_{L^\infty}^2)(\|\nabla \Lambda^s \bd\|^2+\|\Lambda^s\ub\|^2).\no
\end{eqnarray}
\begin{eqnarray}
R_4&\leq& \kappa\|\Lambda^s A\|^2 +C\|\Lambda^s \Grad_\bd W(\bd)\|^2\|\bd\|_{L^\infty}^2\no\\
&\leq& \kappa\|\Lambda^3 A\|^2+C\|\Lambda^s \bd\|^2,\no
\end{eqnarray}
where for $R_4$ we have used Lemma \ref{tr} and \eqref{linfinity}.
The term $R_5, R_6, R_7, R_8$ can be estimated in the same way as $R_1, R_2, R_3, R_4$, respectively.
Besides, the term $R_9, R_{12}$ can be estimated as $R_1$, while $R_{10}, R_{13}$ can be estimated as $R_2$.
Moreover,  from \eqref{critical point of lambda 2} it follows
$$ R_{11}\leq 0.$$
Finally, $K_{11}$ can be estimated similar to $R_4$ so that
\begin{eqnarray}
K_{11}
&\leq& \kappa\|\Delta \Lambda^s \bd\|^2+C\| \Lambda^s \bd\|^2.\no
\end{eqnarray}
Collecting the above estimates together, taking $\kappa>0$ small enough, we conclude
\begin{eqnarray}
&&\frac{d}{dt}(\|\Lambda^s \ub\|^2+\|\nabla \Lambda^s \bd\|^2)+\frac{\mu_4}{2}\|\nabla \Lambda^s \ub\|^2-\frac{1}{\lambda_1}\|\Delta \Lambda^s \bd\|^2
\no\\
&\leq& C(1+\|\nabla \ub\|_{L^\infty}+\|\nabla \bd\|_{L^\infty}^2)\left(1+\|\Lambda^s\ub\|^2+ (1+\|\bd\|_{H^{s-1}}^2)\|\nabla \Lambda^s \bd\|^2\right), \label{dhigh1}
\end{eqnarray}
for $s\geq 3$ and  $t\in [0, T_1]$. Here  $C$ is a constant depending on $\|\ub_0\|_{H^1}$, $\|\bd_0\|_{H^2}$, $T_1$ and the coefficients of the system.

Denote $$Y_s(t)=\|\Lambda^s \ub(t)\|^2+\|\nabla \Lambda^s \bd(t)\|^2.$$
When $s=3$, using the Sobolev embedding theorem, from \eqref{dhigh1} and \eqref{pes1} we infer
\begin{equation}
\frac{d}{dt}Y_3(t)\leq C(1+\|\Lambda^3\ub\|^2+\|\nabla \Lambda^3 \bd\|^2)Y_3(t)\leq C(1+Y_3(t)^2).\label{dhigh2}
\end{equation}
As a result, there exists a time $T_0\in (0, T_1]$ depending on $Y(0)$ (i.e., $\|\ub_0\|_{H^3}$ and $\|\bd_0\|_{H^4}$) such that $Y_3(t)<+\infty$.
Moreover, we get
\begin{equation}
\sup_{0\leq t\leq T_0}\|\ub(t)\|_{H^3}+\|\bd(t)\|_{H^4}\leq C,\no
\end{equation}
which, together with an application of the Sobolev embedding theorem, implies
\begin{equation}
\sup_{0\leq t\leq T_0} \|\nabla \ub\|_{L^\infty}+\|\nabla \bd(t)\|_{L^\infty}.\label{LL}
\end{equation}
Then we can use a simple induction argument. We have shown that $Y_3(t)<+\infty$ on $[0,T_0]$.
Suppose, for $k\geq 3$, we have $Y_k(t)<+\infty$ on $[0,T_0]$. Then  from \eqref{dhigh1} and \eqref{LL} we get
\begin{eqnarray}
\frac{d}{dt}Y_{k+1}(t)&\leq& C(1+\|\nabla \ub\|_{L^\infty}+\|\nabla \bd\|_{L^\infty}^2)\left(1+\|\Lambda^{k+1} \ub\|^2+ (1+\|\bd\|_{H^{k}}^2)\|\nabla \Lambda^{k+1}\bd\|^2\right)\no\\
&\leq& C(1+Y_{k+1}(t)),
\end{eqnarray}
 which yields
 \begin{equation}
 Y_{k+1}(t)\leq (Y_{k+1}(t)+1)e^{Ct}, \quad t\in [0,T_0].\no
 \end{equation}
Thus, for any integer $s\geq 3$, $\ub_0\in H^s_p(Q)$ with $\nabla \cdot \ub_0=0$ and $\bd_0\in H^{s+1}_p(Q)$,
 we have
 \begin{equation}
\sup_{0\leq t\leq T_0}\|\ub(t)\|_{H^s}+\|\bd(t)\|_{H^{s+1}}\leq C,\no
\end{equation}
where $C$ depends on $\|\ub_0\|_{H^s}$, $\|\bd_0\|_{H^{s+1}}$ and the coefficients of the system, while the existence time $T_0>0$ only depends
on $\|\ub_0\|_{H^3}$, $\|\bd_0\|_{H^4}$.
As a consequence, we obtain \textit{a priori} estimates which are uniform
for the approximate solutions on $[0, T_0]$.
By a standard compactness argument, we obtain the existence of the local classical solution $\ub$ on $[0,T_0]$ such that
\begin{equation}
\begin{cases}
& \ub \in L^\infty(0, T_0; H^s_p(Q))\cap L^2(0,
T_0; H^{s+1}_p(Q)),\\
& \bd \in L^\infty(0, T_0; H^{s+1}_p(Q))\cap
L^2(0, T_0; H^{s+2}_p(Q)).\label{rega}
\end{cases}
\end{equation}

The uniqueness of the solution can be easily derived by using the Gronwall inequality in the norm in $L^2\times H^1$, since we have gained
enough higher-order estimates of the solution. Concerning the time continuity property of $\ub, \bd$, we can argue as in \cite{WZ12}.
For any $\delta>0$, we can take $N=N(\delta)$ such that
\begin{equation}
\sup_{0\leq t\leq T_0}\sum_{j>N} 2^{2js} \|\Delta_j \ub(t)\|^2\leq \frac{\delta}{4},\no
\end{equation}
where $\Delta_j$ is the Littlewood--Paley decomposition operator (see, e.g., \cite[Appendix]{WZ12}).
Then, for any $t\in (0,T_0)$ and $\sigma$ such that $t+\sigma\in [0, T_0]$, we have
\begin{eqnarray}
&&\|\ub(t+\sigma)-\ub(t)\|_{H^s}^2\no\\
&\leq& \sum_{j=-1}^{N} 2^{2js} \|\Delta_j \ub(t+\sigma)-\Delta_j \ub(t)\|^2 +2\sup_{0\leq t\leq T_0}\sum_{j>N} 2^{2js} \|\Delta_j \ub(t)\|^2\no\\
&\leq& |\sigma| \sum_{j=-1}^{N} 2^{2js}  \int_0^{T_0}\|\partial_t \ub(t)\|^2 dt+\frac{\delta}{2}\no\\
&\leq& |\sigma|(N+1)2^{2NsN}\int_0^{T_0}\|\partial_t \ub(t)\|^2 dt+\frac{\delta}{2}.\no
\end{eqnarray}
We easily infer from \eqref{rega} and the equation \eqref{e1} that $\ub\in L^2(0, T_0; L^2_p(Q))$. As a result, for $\sigma$ small enough, we have
$\|\ub(t+\sigma)-\ub(t)\|_{H^s}^2\leq \delta$. The time continuity of $\bd$ can be shown in a similar way. The proof is complete.
\end{proof}

\beos
It seems that we are able to prove existence of classical solutions only under the assumptions of \textbf{Case 1}.
This is due to the fact that the Parodi's relation \eqref{lam2} plays an essential role in obtaining differential inequalities for higher-order
Sobolev norms of the solution (see, e.g., \eqref{dhigh1}).
The Parodi's relation leads to important cancelations between higher-order terms like $P_1,...,P_{10}$ and thus the terms with highest-order
derivative cancel so that we can perform the commutator estimates. Indeed, the Parodi's
relation is also crucial to derive the inequality \eqref{A} for the quantity $\mathcal{A}(t)$ (cf. \cite{WXL12} for details).

\eddos


\subsection{Blow-up criterion}

\bete
 \label{blow} {\rm (BKM type blow-up criterion).}
 Suppose that the assumptions of Theorem \ref{locsol} are satisfied. Let $(\ub, \bd)$ be the local classical solution to the E--L system \eqref{e1}--\eqref{e3}
 subject to periodic boundary conditions. Let $T^*$ be the maximal existence time of the solution. If $T^*<+\infty$, then
 \begin{equation}
 \int_0^{T^*}\left( \|\curlu(t)\|_{L^\infty}+\|\nabla \bd(t) \|_{L^\infty}^4 \right) dt=+\infty.   \label{blcon}
 \end{equation}
\ente
In order to prove Theorem \ref{blow}, we need the following lemma.
\bele
 \label{blowlow}
Suppose that the conditions of \textbf{Case 1} are satisfied.
For any $\ub_0\in H^3_p(Q)$ with $\nabla \cdot \ub_0=0$ and $\bd_0\in H^4_p(Q)$, $M>0$ and $T_0>0$, let $(\ub, \bd)$ be a local classical solution to the
E--L system \eqref{e1}--\eqref{e3} subject to periodic boundary conditions. If the following condition is satisfied:
\begin{equation}
\int_0^{T_0}\left( \|\curlu(t)\|_{L^\infty}+\|\nabla \bd(t) \|_{L^\infty}^4 \right) dt\leq M, \label{blcona}
\end{equation}
then
\begin{equation}
\sup_{0\leq t\leq T_0}(\|\ub(t)\|_{H^1}+\|\bd(t)\|_{H^2})\leq C,
\label{u1h2}
\end{equation}
where $C$ is a constant depending on $\|\ub_0\|_{H^1}$, $\|\bd_0\|_{H^2}$, $M$ and the coefficients of the system.
\enle
 \begin{proof}
 In order to obtain the estimate \eqref{u1h2}, we shall derive some energy differential inequality like \eqref{dhigh1}.
 Observe that under assumptions in \textbf{Case 1}, uniform estimates \eqref{eapr1}--\eqref{dAdes} still hold.
 However, here the situation is more delicate because these estimates are rather weak and, in particular, we lose the
 control of $\|\bd\|_{L^\infty(0,T; L^\infty)}$.

Then, applying the curl operator $\nabla \times$ to \eqref{e1}, we obtain
 \begin{equation}
 (\curlu)_t+(\ub \cdot \nabla) \curlu=(\curlu)\cdot \nabla \ub-\nabla \times(\nabla \bd\Delta \bd)+\nabla \times (\nabla \cdot \bsig),\label{cuu}
 \end{equation}
 where we have used the well-known identity
 $\nabla \cdot(\nabla \bd \odot \nabla \bd)=\nabla (\frac12|\nabla \bd|^2)+\nabla \bd\Delta \bd$ and the fact that $\grad \times(\grad\cdot)$ is the null operator.
 Besides, we recall the identity
 \begin{equation}
 \nabla \times \curlu=\nabla (\nabla \cdot \ub)-\Delta \ub,\no
 \end{equation}
 which, together with the incompressibility condition \eqref{e2}, implies
  $$\nabla \times \curlu =-\Delta \ub.$$
 Testing \eqref{cuu} by $\curlu$, we obtain
 \begin{eqnarray}
 && \frac12\frac{d}{dt}\|\curlu\|^2+\frac{\mu_4}{2}\|\Delta \ub \|^2\no\\
 &=& ((\curlu)\cdot \nabla \ub,  \curlu )
  +(\nabla \bd\Delta \bd, \Delta \ub)
    -\mu_1 (\nabla \cdot [(\bd^T A \bd)\bd\otimes \bd], \Delta \ub )\no\\
 &&-\mu_2(\nabla \cdot (\mathcal{N}\otimes \bd),  \Delta \ub)
   -\mu_3 (\nabla \cdot ( \bd\otimes \mathcal{N}),  \Delta \ub)
   \no\\
 &&-\mu_5 (\nabla \cdot[(A\bd)\otimes \bd], \Delta \ub)
 - \mu_6 (\nabla \cdot [\bd\otimes (A \bd)], \Delta \ub)\no\\
 &:=& \sum_{m=1}^7 I_m,
 \no
 \end{eqnarray}
 where we have used the fact
 $$ ((\ub\cdot\nabla )\curlu, \curlu)=\frac12\int_Q(\ub\cdot \nabla )|\curlu|^2dx=0.$$
 Next, applying $\Delta$ to \eqref{e3} and testing the resultant by $\Delta \bd$, we get
 \begin{eqnarray}
&& \frac12\frac{d}{dt} \|\Delta \bd\|^2-\frac{1}{\lambda_1}\|\nabla \Delta \bd\|^2\no\\
&=& -((\Delta \ub \cdot \nabla )\bd,  \Delta \bd) -2\int_Q \nabla_i u_j\nabla_i\nabla_j d_k\Delta d_k dx\no\\
&& +(\Delta (\omega\bd),  \Delta \bd) -\frac{\lambda_2}{\lambda_1}(\Delta (A\bd), \Delta \bd)
 +\frac{1}{\lambda_1}( \Delta \nabla_\bd W(\bd), \Delta \bd)
 \no\\
&:=&\sum_{m=8}^{12} I_m.\no
 \end{eqnarray}
 Obviously, we have the cancellation $I_2+I_8=0$.
 Since the Riesz operators are
bounded in $L^2$ and $\nabla \ub=(-\Delta)^{-1}\nabla (\nabla \times \curlu)$, we deduce $\|\nabla \ub\|\leq C\|\curlu\|$.
Then for $I_1$ we see that
 $$
 I_1\leq  \|\curlu\|_{L^\infty}\|\nabla \ub\|\|\curlu\|\leq  C\|\curlu\|_{L^\infty}\|\curlu\|^2.
 $$
Moreover, using the Gagliardo--Nirenberg inequality and the H\"older inequality, for $I_9$ and $I_{12}$ we have
 \begin{eqnarray}
 I_9&\leq& C\|\nabla \ub\|\|\nabla ^2\bd\|_{L^4}^2\no\\
 &\leq& C\|\curlu\|(\|\nabla \bd\|_{L^\infty}\|\nabla\Delta \bd\|+\|\nabla \bd\|_{L^\infty}^2)\no\\
&\leq& \kappa \|\nabla \Delta\bd\|^2+ C\|\nabla \bd\|_{L^\infty}^2\|\curlu\|^2, \no
 \end{eqnarray}
 \begin{eqnarray}
 I_{12}
 &=& \frac{1}{\varepsilon}\int_Q \nabla  (|\bd|^2\bd) : \nabla \Delta \bd dx+ \frac{1}{\varepsilon}\|\Delta \bd\|^2\no\\
 &\leq& C \|\bd\|_{L^6}^2\|\nabla \bd\|_{L^6}\|\nabla \Delta \bd\|+\frac{1}{\varepsilon}\|\Delta \bd\|^2\no\\
 &\leq& \kappa \|\nabla\Delta \bd\|^2+ C\|\bd\|_{H^1}^4\|\nabla^2 \bd\|^2+\frac{1}{\varepsilon}\|\Delta \bd\|^2\no\\
 &\leq& \kappa  \|\nabla\Delta \bd\|^2+ C\|\Delta \bd\|^2.\no
\no
 \end{eqnarray}
By integration by parts and the fact that $\omega$ is antisymmetric, from the H\"older inequality and the  Young inequality we deduce
 \begin{eqnarray}
 I_3&=& \mu_1([(\bd^T A\bd)\bd\otimes\bd], \Delta \nabla \ub)\no\\
  &=&  \mu_1([(\bd^T A\bd)\bd\otimes\bd], \Delta (A+\omega))
  \no\\
  &=& \mu_1([(\bd^T A\bd)\bd\otimes\bd], \Delta A)\no\\
  &=& -\mu_1\int_Q \nabla_l(d_kA_{kp}d_p d_id_j)\nabla_l A_{ij} dx\no\\
  &=&-\mu_1\int_{Q}(d_kd_p\nabla_lA_{kp})^2dx\no\\
  &&-\mu_1\int_{Q}A_{kp}(\nabla_ld_kd_p +d_k\nabla d_p)d_id_j \nabla_lA_{ij}dx
\no\\
&&-\mu_1\int_{Q}A_{kp}d_kd_p( d_i\nabla_ld_j+ \nabla_l d_i d_j)\nabla_lA_{ij}dx\no\\
  &\leq& -\frac{\mu_1}{2} \int_{Q}(d_kd_p\nabla_lA_{kp})^2dx+\kappa \|\Delta \ub\|^2\no\\
      &&
      + C\mu_1(\|\bd\|_{L^\infty}^2+\|\bd\|_{L^\infty}^6)\|\nabla \ub\|^2\|\nabla  \bd\|_{L^\infty}^2.\no
  \end{eqnarray}
 Next, consider the following inequality
 $$
 \|\bd\|_{L^\infty}\leq C\|\nabla \bd\|_{L^\infty}^\frac13\|\bd\|_{L^6}^\frac23+C\|\bd\|_{L^6}.
 $$
Then we infer
  \begin{eqnarray}
  && (\|\bd\|_{L^\infty}^2+\|\bd\|_{L^\infty}^6)\|\nabla \ub\|^2\|\nabla \bd\|_{L^\infty}^2\no\\
  &\leq& C\left[1+\Big(\|\nabla \bd\|_{L^\infty}^\frac13\|\bd\|_{L^6}^\frac23+C\|\bd\|_{L^6}\Big)^6\right]\|\curlu\|^2\|\nabla \bd\|_{L^\infty}^2\no\\
  &\leq& C\Big(1+\|\nabla \bd\|_{L^\infty}^4\Big)\|\curlu\|^2.\no
  \end{eqnarray}
  As a consequence, we have
  \begin{eqnarray}
  I_3&\leq& -\frac{\mu_1}{2} \int_{Q}(d_kd_p\nabla_lA_{kp})^2dx
  +\kappa \|\Delta \ub\|^2 + C\mu_1\Big(1+\|\nabla \bd\|_{L^\infty}^4\Big)\|\curlu\|^2.
 \no
  \end{eqnarray}
 Concerning $I_4, ..., I_7$, using integration by parts, we obtain (cf. e.g., \cite[Appendix]{WXL12})
   \begin{equation}
   I_4+I_5=(\mu_2+\mu_3)\int_Q d_j\mathcal{N}_i\Delta A_{ij}dx+(\mu_2-\mu_3)(\mathcal{N},\Delta
\omega \,\bd)
 \label{expansion of mu 2 and 3}
 \end{equation}
 and
 \begin{eqnarray}
  I_6+I_7
&=&-(\mu_5+\mu_6)\int_{Q}|d_j\nabla_lA_{ji}|^2dx-(\mu_5+\mu_6)\int_{Q}\nabla_ld_jd_kA_{ki}
\nabla_lA_{ij}dx
\no\\
&&-(\mu_5+\mu_6)\int_{Q}d_j\nabla_ld_kA_{ki}
\nabla_lA_{ij}dx+(\mu_5-\mu_6)\big(A\bd, \Delta\omega \bd \big)\no\\
&:=& J_1+J_2+J_3+(\mu_5-\mu_6)\big(A\bd, \Delta\omega \bd \big).
\label{expansion of mu 5 and mu 6}
 \end{eqnarray}
On account of equation \eqref{e3} and integrating by parts, we get
 \begin{eqnarray}
  I_{10}
 &=&(\Delta \bd, \Delta\omega\,
\bd)+2\int_Q\Delta d_i \nabla_l\omega_{ij}\nabla_l d_j
dx+(\Delta \bd,
\omega\Delta \bd)
   \no\\
&=&-\lambda_1\int_{Q}d_j \mathcal{N}_i \Delta\omega_{ij}\, dx-\lambda_2\big(A\bd,
\Delta\omega \bd \big)+  \big(\nabla_{\bd}W(\bd),
\Delta\omega \bd \big)
 \no\\
&&  +2\int_Q\Delta d_i
\nabla_l\omega_{ij}\nabla_l d_j dx+(\Delta\bd, \omega\Delta \bd)
  \no\\
&=& -\lambda_1(\mathcal{N}, \Delta\omega\,\bd)-\lambda_2\big(A\bd, \Delta\omega \bd
\big) -\int_Q\nabla_l[ (\nabla_{\bd}W(\bd))_i d_j]\nabla_l\omega_{ij}dx  \no\\
&&  -\int_Q\nabla_l \Delta d_i \omega_{ij}\nabla_l d_j dx
 +\int_Q\Delta d_i \nabla_l\omega_{ij}\nabla_l d_j dx\no\\
 &:=& -\lambda_1(\mathcal{N}, \Delta\omega\,\bd)-\lambda_2\big(A\bd, \Delta\omega \bd
\big)+J_4+J_5+J_6,
 \label{expansion of laplace omega d}
 \end{eqnarray}
 and
 \begin{eqnarray}
 I_{11}
&=&\lambda_2 \big(\mathcal{N}, \Delta(A\bd)
\big)+\frac{(\lambda_2)^2}{\lambda_1}\big(A\bd,
\Delta(A\bd) \big)-\frac{\lambda_2}{\lambda_1}(\nabla_\bd W(\bd),\Delta(A\bd) \big))
  \no\\
&=&\lambda_2\int_Q \mathcal{N}_i\Delta A_{ij}d_jdx+2\lambda_2\int_Q
\mathcal{N}_i\nabla_l A_{ij}\nabla_l d_j dx+\lambda_2(\mathcal{N}, A\Delta \bd)
  \no\\
&& -\frac{(\lambda_2)^2}{\lambda_1}\int_{Q}|\nabla_l
(A_{ij}d_j)|^2dx+\frac{\lambda_2}{\lambda_1}(\nabla  \nabla_\bd W(\bd),\nabla (A\bd) \big)\no\\
&=& \lambda_2\int_Q \mathcal{N}_i\Delta A_{ij}d_jdx
 + \lambda_2\int_Q \mathcal{N}_i\nabla_l A_{ij}\nabla_l d_j
dx -\lambda_2\int_Q \nabla_l\mathcal{N}_i A_{ij}\nabla_l d_j dx\no\\
&&  -\frac{(\lambda_2)^2}{\lambda_1}\int_{Q}|d_j\nabla_lA_{ji}|^2dx
-\frac{(\lambda_2)^2}{\lambda_1}\int_{Q}|A_{ij}\nabla_l d_j|^2dx
\no\\
&& -\frac{2(\lambda_2)^2}{\lambda_1}\int_{Q}\nabla_l A_{ij} d_j
A_{ik}\nabla_l d_k
dx +\frac{\lambda_2}{\lambda_1}  \int_Q \nabla_l (\nabla_\bd W(\bd))_i\nabla_lA_{ij}d_j dx\no\\
&& + \frac{\lambda_2}{\lambda_1}  \int_Q \nabla_l (\nabla_\bd W(\bd))_i A_{ij} \nabla_l d_j dx
\no\\
&=& \lambda_2\int_Q \mathcal{N}_i\Delta A_{ij}d_jdx -\frac{(\lambda_2)^2}{\lambda_1}\int_{Q}|d_j\nabla_lA_{ji}|^2dx
\no\\
&&
+\frac{\lambda_2}{\lambda_1}\int_Q
\nabla_l \Delta d_iA_{ij}\nabla_l d_j dx-\frac{\lambda_2}{\lambda_1}\int_Q \Delta d_i\nabla_l
A_{ij}\nabla_l d_j dx\no\\
&& -\frac{(\lambda_2)^2}{\lambda_1}\int_{Q}\nabla_l A_{ij} d_j
A_{ik}\nabla_l d_k
dx
 +\frac{\lambda_2}{\lambda_1}  \int_Q \nabla_l [(\nabla_\bd W(\bd))_i d_j]\nabla_lA_{ij} dx \no\\
\no\\
&:=& \lambda_2\int_Q \mathcal{N}_i\Delta A_{ij}d_jdx + J_7+...+J_{11}.
 \label{expansion of laplace A d}
 \end{eqnarray}
By the Parodi's relation
\eqref{lam2} and the condition \eqref{lama1}, the first term of the right-hand side of
\eqref{expansion of laplace A d} cancels with the first term of the
right-hand side of \eqref{expansion of mu 2 and 3}. Besides, due to the relation  \eqref{lama1}, the first term on the right-hand side of
\eqref{expansion of laplace omega d} cancels with the second term of
the right-hand side of \eqref{expansion of mu 2 and 3} and the
second term on the right-hand side of \eqref{expansion of laplace
omega d} cancels with the fourth term of the right-hand side of
\eqref{expansion of mu 5 and mu 6}.

 Keeping in mind these special cancellations between the nonlinear terms of the highest-order derivatives, we are able to estimate the remaining
 terms in $I_6, I_7, I_{10}, I_{11}$, namely, $J_1, ..., J_{12}$.

 It follows from \eqref{mu56} and \eqref{critical point
of lambda 2} that $$J_1+J_7\leq 0.$$
 Next, we have
 \begin{eqnarray}
 J_2+J_3&\leq& C(\mu_5+\mu_6)\|\Delta \ub\|\|\nabla \ub\|\|\nabla \bd\|_{L^\infty}\|\bd\|_{L^\infty}\no\\
 &\leq& \kappa\|\Delta \ub\|^2+ C(\mu_5+\mu_6)^2 \Big(\|\nabla \bd\|_{L^\infty}^\frac13\|\bd\|_{L^6}^\frac23+\|\bd\|_{L^6}\Big)^2 \|\nabla \bd\|_{L^\infty}^2\|\curlu\|^2\no\\
 &\leq& \kappa\|\Delta \ub\|^2+C(\mu_5+\mu_6)^2\Big(1+\|\nabla \bd\|_{L^\infty}^\frac83\Big)\|\curlu\|^2\no,
 \end{eqnarray}
  \begin{eqnarray}
 J_4&\leq& C \|\Delta \ub\|\| \nabla \bd\|_{L^6}(\|\bd\|_{L^3}+ \|\bd\|_{L^9}^3) \no\\
 &\leq& \kappa \|\Delta \ub\|^2+ C\left(1+\|\nabla \bd\|_{L^\infty}^\frac13\right)\| \Delta \bd\|^2.\no
 \end{eqnarray}
 For $J_5, J_6$, by an integration by parts and the fact that $\omega$ is anti-symmetric, we get
 \begin{eqnarray}
 J_5+J_6&=& 2\int_Q\Delta d_i \nabla_l\omega_{ij}\nabla_l d_j dx \leq \kappa\|\Delta \ub\|^2+ C\|\nabla \bd\|_{L^\infty}^2\|\Delta \bd\|^2.\no
 \end{eqnarray}
  Integrating by parts and using the H\"older inequality and the interpolation inequalities
  \begin{eqnarray}
 &&  \|\Delta \bd\|_{L^4}\leq C( \|\nabla \Delta \bd\|^\frac12\|\nabla \bd\|_{L^\infty}^\frac12+\|\nabla \bd\|_{L^\infty}),
 \no\\
 && \|\nabla \bd\|_{L^\infty}\leq C\|\nabla \bd\|_{H^1}^\frac12\|\nabla \bd\|_{H^2}^\frac12,
 \no
  \end{eqnarray}
 then  we get
\begin{eqnarray}
 J_8+J_9
&=& -\frac{\lambda_2}{\lambda_1}\int_Q \Delta d_i
A_{ij}\Delta d_j dx
-\frac{2\lambda_2}{\lambda_1}\int_Q
 \Delta d_i \nabla_l A_{ij}\nabla_l d_j dx\no\\
 &\leq& |\lambda_2|\|\nabla \ub\|\|\Delta \bd\|_{L^4}^2+|\lambda_2| \|\Delta \ub\|\|\nabla \bd\|_{L^\infty}\|\Delta \bd\|\no\\
&\leq& \kappa\|\Delta \ub\|^2+\kappa \|\nabla \Delta \bd\|^2+C(\lambda_2)^2 \|\nabla \bd\|_{L^\infty}^2(\|\curlu\|^2+\|\Delta \bd\|^2)\no\\
&&  +C\|\nabla  \bd\|^2_{L^\infty}.\no
\end{eqnarray}
For $J_{10}$, we have
 \begin{eqnarray}
 J_{10}&\leq & \frac{(\lambda_2)^2}{\lambda_1}\|\Delta \ub\|\|\bd\|_{L^\infty}\|\nabla \ub\|\|\nabla \bd\|_{L^\infty}\no\\
 &\leq& \kappa\|\Delta \ub
\|^2+C (\lambda_2)^4\|\bd\|_{L^\infty}^2\|\nabla \bd\|^2_{L^\infty}\|\curlu\|^2\no\\
&\leq& \kappa\|\Delta \ub
\|^2+C (\lambda_2)^4(1+\|\nabla \bd\|^\frac23_{L^\infty})\|\nabla \bd\|^2_{L^\infty}\|\curlu\|^2.\no
  \end{eqnarray}
Finally, $J_{11}$ can be estimated just as $J_4$.

Collecting the above estimates and taking $\kappa>0$ small enough,  from the the Young inequality and \eqref{critical point of lambda 2}
(i.e., $|\lambda_2|$ can be bounded by $\mu_5+\mu_6\geq 0$), we infer
 \begin{equation}
\frac{d}{dt}\left(\|\curlu\|^2+\|\Delta \bd\|^2\right)
 \leq G(t)(1+\|\curlu\|^2+\|\Delta \bd\|^2),
 \label{diff}
\end{equation}
where
 \begin{equation}
 G(t)=C\left[1+\|\curlu\|_{L^\infty}+\|\nabla \bd\|_{L^\infty}^2
 +(\mu_5+\mu_6)^2\|\nabla \bd\|_{L^\infty}^\frac83+\mu_1\|\nabla \bd\|_{L^\infty}^4\right].\label{GG}
 \end{equation}
Here $C$ is a constant depending on $\|\ub_0\|$, $\|\bd_0\|_{H^1}$ and the coefficients of the system. Finally, by the Gronwall inequality,
we deduce
 \begin{eqnarray}
 \|\curlu(t)\|^2+ \|\Delta \bd(t)\|^2\leq  \left(1+\|\curlu_0\|^2+ \|\Delta \bd_0\|^2\right)
 {\rm exp}\left(\int_0^t G(\tau) d\tau\right), \quad t\in [0,T_0].\no
 \end{eqnarray}
 This and the lower order estimates \eqref{eapr1}--\eqref{eapr2} yield the conclusion \eqref{u1h2}. The proof is complete.
 \end{proof}

\par \textbf{Proof of Theorem \ref{blow}}. We prove the theorem by a contradiction argument. Suppose that \eqref{blcon} is not true.
Then there exists a positive constant $ M$ such that
 \begin{equation}
\int_0^{T^*}\left( \|\curlu(t)\|_{L^\infty}+\|\nabla \bd(t) \|_{L^\infty}^4 \right) dt\leq M,\label{blbd}
 \end{equation}
 which yields (cf. Lemma \ref{blowlow})
 \begin{equation}
\sup_{0\leq t\leq T^*}(\|\ub(t)\|_{H^1}+\|\bd(t)\|_{H^2})\leq C, \label{u1h3}
\end{equation}
where $C$ is a constant depending on $\|\ub_0\|_{H^1}$, $\|\bd_0\|_{H^2}$, $M$ and the coefficients of the system.

Recalling the proof of Theorem \ref{locsol} and using the bound \eqref{u1h3}, we are able to derive the following differential inequality \eqref{dhigh1} for the quantity $Y_3(t)$
 \begin{equation}
 \frac{d}{dt} Y_3(t)\leq C \left(1+\|\nabla  \ub\|_{L^\infty}+\|\nabla \bd\|_{L^\infty}^2\right)(Y_3(t)+1).\label{YYY}
 \end{equation}
Combining the critical logarithmic Sobolev inequality (cf. \cite{BKM})
 $$
 \|\nabla \ub\|_{L^\infty}\leq C\left(1+\|\curlu\|+\|\curlu\|_{L^\infty} \ln(e+\|\ub\|_{H^3})\right),
 $$
 with estimate \eqref{u1h3}, we get
 $$
 \|\nabla \ub\|_{L^\infty}\leq C\left(1+\|\curlu\|_{L^\infty}\ln(e+\|\ub\|_{H^3})\right).
 $$
 Then, from \eqref{YYY} we infer
 \begin{eqnarray}
 \frac{d}{dt}(Y_3(t)+e+1)\leq C(1+\|\curlu\|_{L^\infty}+\|\nabla \bd\|_{L^\infty}^2)(Y_3(t)+e+1)\ln(Y_3(t)+e+1).\no
 \end{eqnarray}
By using the Gronwall inequality, for $t\in [0,T^*)$, we conclude
\begin{equation}
Y_3(t)+e+1\leq {\rm exp} \left( \ln[Y_3(0)+e+1]{\rm exp}\left(C\int_0^t (1+\|\curlu\|_{L^\infty}+\|\nabla \bd\|_{L^\infty}^2)d\tau\right)\right).\no
\end{equation}
 Therefore, if $T^*<+\infty$ and \eqref{blbd} holds, an application of the Young inequality gives
 $$\int_0^{T^*} (1+\|\curlu\|_{L^\infty}+\|\nabla \bd\|_{L^\infty}^2)dt\leq
 \int_0^{T^*} (1+\|\curlu\|_{L^\infty}+\|\nabla \bd\|_{L^\infty}^4)dt<+\infty.$$
 Then we deduce that  $Y_3(t)$ is bounded for $t\in [0,T^*)$ and the local classical solution $(\ub, \bd)$ can be extended beyond $t=T^*$.
 This leads to a contradiction of the definition of maximal existence time $T^*$. The proof is complete. \qed

 \beos
Comparing our blow-up criterion \eqref{blcon} to the one obtained in \cite{HW12} for the simplified liquid crystal system with constraint $|\bd|=1$ such that
 \begin{equation}
\int_0^{T^*} (\|\curlu(t)\|_{L^\infty}+\|\nabla \bd(t)\|_{L^\infty}^2) dt=+\infty,\label{blaa}
\end{equation}
they differ only on the power of $\|\nabla \bd\|_{L^\infty}$. Observe that in \cite{HW12} the term $\|\nabla \bd\|_{L^\infty}$
has to be included in the criterion mainly because of the difficulty due to heat flow of harmonic maps.
For the full E--L system with constraint $|\bd|=1$ (cf. \cite{WZZ12}), higher-order stress tensors will not introduce further troubles and the same criterion \eqref{blaa}
still works. However, in our case, the power of $\|\nabla \bd\|_{L^\infty}$ has to be larger than those in \cite{WZZ12, HW12}, due to the lack of control on $\|\bd\|_{L^\infty}$.
Furthermore, the expression \eqref{GG} indicates that for different choices of the Leslie coefficients (always under assumptions of \textbf{Case 1}), the blow-up criterion \eqref{blcon}
can be improved.
For instance, under the special case $\mu_1= \mu_5+\mu_6=0$ (e.g., taking $\alpha=\frac12$ in the simplified liquid crystal system studied in \cite{CR, ffrs, wuxuliu}),
\eqref{blcon} can be replaced by \eqref{blaa}. In this case, due to \eqref{critical point of lambda 2}, we have also $\lambda_2=0$, namely, the stretching effect is neglected
in the E--L system while the rotation effect is kept. When $\alpha=1$, we refer to \cite{ZF10} for another type of criterion only in terms of the velocity
 \begin{equation}
 \int_0^{T^*} \frac{\|\nabla \ub(t)\|_{L^p}^r}{ 1+ \ln(e+\|\nabla \ub(t)\|_{L^p})}dt<+\infty, \quad \frac{2}{r}+\frac{3}{p}= 2, \ 2\leq p\leq 3,\no
 \end{equation}
 but with a restricted range of $p$ instead of the natural one $\frac32\leq p\leq \infty$.

If we neglect the parallel transport of the director $\bd$ (i.e., the term $-\omega\bd+\frac{\lambda_2}{\lambda_1}A\bd$ in \eqref{e3}) and drop all the other Leslie stress terms in
$\bsig$ except the fluid viscosity term $\mu_4 A $, then the E--L system \eqref{e1}--\eqref{e3} will be reduced to the simplest model studied in \cite{LinLiusimply}.
For this case, a logarithmical improved BKM's criterion in terms of the velocity field $\ub$ can be obtained (cf. \cite{LZC12}), meaning
$$
\int_0^{T^*} \frac{\|\curlu(t)\|_{BMO}}{\sqrt{\ln(e+\|\curlu(t)\|_{BMO})}}dt<+\infty.
$$
We recall here the definition of BMO (Bounded Mean Oscillation) space
$$ {\rm BMO} = \left\{ f \in L_{loc}(\mathbb{R}^3) \, : \, \|f\|_{BMO} = \sup_{R>0,\, x \in \mathbb{R}^3} \frac{1}{|B_R(x)|} \int_{B_R(x)}|f(y) - \bar f_{B_R(x)}|dy < \infty \right \},\mathcal{\mathbb{}}$$
where $\bar f_{B_R(x)}$ stands for the average of $f$ over $B_R(x)$.
 \eddos

\medskip

\noindent \textbf{Acknowledgements.} The authors are very grateful
to the referee for her/his helpful comments and suggestions.


\end{document}